\documentclass[11pt]{amsart}
\usepackage{amsmath,amsthm, amssymb, amsfonts,amscd, epsfig}

\usepackage[mathscr]{eucal}

\theoremstyle{plain}

\newcommand{\rank}{\mathrm{Rank}\;}
\newcommand{\Span}{\mathrm{Span}\;}
\newtheorem{thm}{Theorem}[section]
\newtheorem{lem}[thm]{Lemma}
\newtheorem{pro}[thm]{Proposition}
\newtheorem{cor}[thm]{Corollary}

\theoremstyle{definition}
\newtheorem{Def}[thm]{Definition}
\theoremstyle{remark}
\newtheorem{rem}[thm]{Remark}
\newtheorem{exa}[thm]{Example}

\newcommand{\sat}{\mathrm{sat}}
\newcommand{\bideg}{\mathrm{bideg}}

\newcommand{\mc}[1]{\mathcal{#1}}

\newcommand{\mbf}[1]{\mathbf{#1}}

\newcommand{\coh}[3]{ H ^{#1}  ( #2, \,  #3 )}

\newcommand{\Kos}{\mathrm{Kos}\;}
\newcommand{\Syz}{\mathrm{Syz}\;}

\newcommand{\Reg}{{\mathrm{Reg}}}

\newcommand{\V}{{\mathbb V}}
\newcommand{\Z}{{\mathbb Z}}

\newcommand{\C}{{\mathbb C}}

\newcommand{\im}[1]{{\rm Im}\,(#1)}
\newcommand{\Ker}[1]{{\rm Ker}\,(#1)}
\newcommand{\proj}[1]{{\mathbf P}^{#1}}

\newcommand{\set}[1]{\left\{#1\right\}}

\begin{document}

\title[Implicitization]{Equations of parametric surfaces with
  base points via syzygies}

\author{William A. Adkins, J. William Hoffman and Hao Hao Wang }
\address{Department of Mathematics\\
              Louisiana State University\\
              Baton Rouge, Louisiana 70803}

\email{adkins@math.lsu.edu, hoffman@math.lsu.edu, wang\_h@math.lsu.edu}

\begin{abstract}
Let $S$ be a parametrized surface in $\proj{3}$ given as the image
of $\phi: \proj{1} \times \proj{1} \to \proj{3}$.  This paper will
show that the use of syzygies in the form of a combination of
moving planes and moving quadrics provides a valid method for
finding the implicit equation of $S$ when certain base points are
present. This work extends the algorithm provided by Cox \cite{C}
for when $\phi$ has no base points, and it is analogous to some of
the results of Bus\'e, Cox, and D'Andrea \cite{BCD} for the case
when $\phi: \proj{2} \to \proj{3}$ has base points.
\end{abstract}

\thanks{We would like to thank David Cox for
numerous discussions and suggestions.}


\maketitle

\section{Introduction} \label{s:intro}

The use of syzygies has been explored in a number of recent works
as an alternative to resultants for producing determinantal
formulas for the equations of rationally parametrized curves and
surfaces.  The  article by Cox \cite{C1} provides a detailed
survey of the current status of the problem of finding the
implicit equation of a rational surface $S\subset \proj{3}$
described implicitly by a map $\phi: X\to \proj{3}$, where $X$ is
either $\proj{2}$ or $\proj{1}\times \proj{1}$.  The reader is
referred to this paper and its references for a discussion of the
history of the use of syzygies in the implicitization problem,
that is, the problem of finding a generator for the ideal $I(S)$
from the knowledge of $\phi$.

In this paper we will only consider the case $X=\proj{1}\times
\proj{1}$, so that our parametrization map $\phi:\proj{1}\times
\proj{1}\to \proj{3}$ will have the form
\begin{equation}
\phi=[a_0(s,\,u;\,t,\,v),\, a_1(s,\,u;\,t,\,v),\,
a_2(s,\,u;\,t,\,v),\,a_3(s,\,u;\,t,\,v)],
\end{equation}
where $a_0$, $a_1$, $a_2$,  $a_3\in R=\C[s,\, u,\, t,\,v]$  are
bihomogeneous polynomials of bidegree $(m,\,n)$.  Moreover, we
always assume that $\gcd(a_0,\, a_1,\,a_2,\,a_3)=1$.  Even with
the gcd assumption, it can happen that there are points of
$\proj{1}\times \proj{1}$ where all of the $a_i$, $0\le i\le 3$,
vanish simultaneously. These are points, referred to as \emph{base
points},  where the map $\phi$ is not defined. The goal of this
paper is to study the implicitization problem when base points are
present, but we will start by summarizing how syzygies were
employed by Cox, Goldman, and Zhang \cite{CGZ}, \cite{C} to
produce a determinantal equation for $S$ in the case of no base
points. If $\phi$ has no base points, that is, $\V(a_0, a_1, a_2,
a_3)=\emptyset$ in $\proj{1} \times \proj{1}$, and $\phi$ is
generically one-to-one, then the image of $\phi$ is a surface
$S\subset \proj{3}$ of degree $2mn$, \cite[Theorem 3.1]{C},
\cite{CGZ}.

In the polynomial ring
$\C[s,u,t,v,x_0,x_1,x_2,x_3]=R[x_0,x_1,x_2,x_3]$, consider the
polynomial $\sum_{i=0}^3A_ix_i$ where $A_i \in R$ $(0\le i\le 3)$
are bihomogeneous polynomials, all of the same bidegree.  If we
fix a point $\mathbf{p}=[s,u;t,v] \in \proj{1} \times \proj{1}$,
then $\sum_{i=0}^3 A_i(\mathbf{p})x_i=0$ is an equation of a plane
in $\proj{3}$, provided some $A_i(\mathbf{p})\ne 0$. When the
point $\mathbf{p}$ changes, we will obtain different equations of
planes in $\proj{3}$. This suggests the following definition:

\begin{Def}
A \textit{moving plane} on $\proj{3}$ is a polynomial of the form
\[\sum_{i=0}^3A_ix_i\]
where, for $0\le i\le 3$, $x_i$ are homogeneous coordinates on
$\proj{3}$ and $A_i \in R$ are bihomogeneous polynomials of the
same bidegree $(k,l)$, which we will call the \emph{bidegree of
the moving plane}. We say the \textit{moving plane follows the
parametrization} $\phi$ if
\[ \sum_{i=0}^3A_i(\mathbf{p})a_i(\mathbf{p})=0,\quad
\text{for all $\mathbf{p} \in \proj{1}\times \proj{1}$},\] which
is equivalent to
\begin{equation}\label{defmp}
\sum_{i=0}^3A_ia_i=0 \in \C[s,u,t,v]
\end{equation}
where the polynomials $a_i$ $(0\le i\le 3)$ are the parameters
that define the surface $S$.
\end{Def}
In the language of commutative algebra,  Equation \eqref{defmp}
states that the moving plane $\sum_{i=0}^3A_ix_i$ follows the
parametrization $\phi$ if and only if
\[ (A_0, A_1, A_2, A_3)\in \Syz(a_0, a_1, a_2, a_3)\]
where $\Syz(a_0, a_1, a_2, a_3)$ denotes the syzygy submodule of
$R^4$ determined by $a_0, a_1, a_2, a_3\in R$.

Analogously:
\begin{Def}
A \textit{moving quadric} is polynomial
\begin{equation}\label{defmq}
 \sum_{0\leq i
\leq j \leq 3} A_{ij}x_ix_j
\end{equation}
which is quadratic in the homogeneous variables $x_i$ ($0\le i\le
3$) and where all of the $A_{ij}\in R$ are bihomogeneous
polynomials of the same bidegree $(k,l)$. We will call this common
bidegree $(k,l)$ the \emph{bidegree of the moving quadric}.

As with moving planes, a \textit{moving quadric follows the
parametrization} $\phi$, if
$$
 (A_{00}, A_{01}, \ldots, A_{33})\in \Syz(a^2_0,a_0a_1,\ldots,
a^2_3),
$$
which means that
$$ \sum_{0\leq i
\leq j \leq 3}
A_{ij}(\mathbf{p})a_i(\mathbf{p})a_j(\mathbf{p}),\quad \text{for
all $\mathbf{p} \in \proj{1}\times \proj{1}$}.$$
\end{Def}
We will primarily have occasion to focus on moving planes and
moving quadrics of bidegree $(m-1,n-1)$ that follow the
parametrization $\phi$, which we have assumed has bidegree
$(m,n)$. If $R_{k,l}\subset R$ denotes the bihomogeneous forms of
bidegree $(k,l)$, then the moving planes of  bidegree $(m-1,n-1)$
make up the kernel of the complex linear map
\[\begin{CD}MP: R_{m-1,n-1}^4 @>{[a_0\;a_1\;a_2\;a_2\;a_3]}>>
R_{2m-1,2n-1}\end{CD}\] given by
\[ MP(A_0, A_1, A_2, A_3) = \sum_{i=0}^3A_ia_i.\]
Note that  the standard monomial basis of $R_{k,l}$ is
$$\mathcal{B}_{k,\,l}=\set{s^iu^{k-i}t^jv^{l-j}: 0\le i\le k,\;
0\le j\le l}$$ so that $\dim_{\C}R_{k,l}=(k+1)(l+1)$. With respect
to the standard bases $\mathcal{B}^4_{m-1,\,n-1}$ on
$R_{m-1,\,n-1}^4 $ and $\mathcal{B}_{2m-1,\,2n-1}$ on
$R_{2m-1,\,2n-1}$,  the linear map $MP$ is represented by a $4mn
\times 4mn$ matrix that, by abuse of notation, we will also denote
by  $MP$. If $\phi$ has no base points and is generically
one-to-one, then $MP$ is an isomorphism \cite[page 8]{C},
 so that there are no moving planes of bidegree
$(m-1,n-1)$.  One of our  results (Lemma \ref{L:mp}) is the
verification that certain base points of total multiplicity $k$
will will have the effect of producing exactly $k$ linearly
independent moving planes.

Similarly, the moving quadrics of  bidegree $(m-1, n-1)$ are the
kernel of the map
\[\begin{CD}MQ: R_{m-1,n-1}^{10} @>{[a^2_0\ a^{ }_0a^{ }_1\
\cdots\  a^2_3]}>> R_{3m-1,3n-1}\end{CD}\] given by
\[ MQ(A_{00},A_{01},\ldots,A_{33}) =
\sum_{0 \leq i \leq j \leq 3} A_{ij}a_ia_j.\] As for the case of
moving planes, we will  identify the map $MQ$ with  the $9mn\times
10mn$ matrix which represents $MQ$ in the standard bases. Since
$$\dim R^{10}_{m-1,n-1}-\dim R_{3m-1,3n-1} =10mn-9mn=mn,$$ it
follows that $\dim \Syz(a_0^2,\cdots, a_3^2)_{m-1,n-1} \ge mn$ and
\[ \dim \Syz(a_0^2,\cdots, a_3^2)_{m-1,n-1} =mn \Longleftrightarrow
MQ \text{ has maximal rank}.\] Thus, if $MQ$ has maximal rank, we
can choose a basis of exactly $mn$ linearly independent moving
quadrics of bidegree $(m-1,n-1)$ which follow the parametrization
$\phi$. Each of these $mn$ linearly independent moving quadrics
$Q_k$ $(1 \le k \leq mn)$ can be written as
\begin{eqnarray}
Q_k &=&  \sum_{0\leq i \leq j \leq 3} A_{ij} x_ix_j\notag\\
&=&\sum_{0\leq i \leq j \leq 3} \left(\sum_{\alpha=0}^{m-1}
\sum_{\beta=0}^{n-1}A_{ij,\alpha \beta}s^{\alpha}t^{\beta}
\right)x_ix_j\notag\\
&=& \sum_{\alpha=0}^{m-1}\sum_{\beta=0}^{n-1}\left(\sum_{0\leq i
\leq j
\leq 3} A_{ij, \alpha \beta}x_ix_j\right)s^{\alpha}t^{\beta}\notag \\
\label{E:quad}&=&
\sum_{\alpha=0}^{m-1}\sum_{\beta=0}^{n-1}Q_{k,\alpha \beta}
(x_0,x_1,x_2,x_3)s^{\alpha}t^{\beta}
\end{eqnarray} where $Q_{k,\alpha \beta}$ is a quadric in $x_i$ with
coefficients in $\C$.  To simplify the notation somewhat, we have
identified the bihomogeneous monomial
$s^{\alpha}u^{m-1-\alpha}t^{\beta}v^{n-1-\beta}$ with its
particular dehomogenized form $s^{\alpha }t^{\beta}$ obtained by
taking $u=v=1$. Arrange the $Q_{k, \alpha \beta}$ into a square
matrix $M$ of size $mn\times mn$, where the columns of the matrix
$M$ are indexed by the monomial basis
$\{s^{\alpha}t^{\beta}\}_{\alpha=0, \beta=0}^{m-1, n-1}$ of
$R_{m-1,n-1}$, and the rows are indexed by the $mn$ moving
quadrics $Q_k$ ($1\leq k \leq mn$).  Since each entry of $M$ is a
quadric in $x_i$,  we may write
\[ M=[Q_{k, \alpha \beta}], \]
so that the determinant of $M$, denoted as usual by $|M|$, is a
polynomial in the variables $x_i$ of degree $\leq 2mn$. One of the
main results of \cite{CGZ} uses the matrix $M$ to give a
determinantal equation for $S=\im{\phi}$.
\begin{thm}
Suppose that $\phi: \proj{1}\times \proj{1} \to \proj{3}$ has no
base points and is generically one-to-one. If $MP$ has maximal
rank, then so does $MQ$ and furthermore, the image of $\phi$ is
defined by the determinantal equation $|M|=0$.
\end{thm}
\begin{proof}
See \cite[Theorem 3.1]{C}.
\end{proof}
The goal of this paper is to prove a similar result where base
points are allowed so long as each base point is a local complete
intersection  and the total multiplicity of all base points does
not exceed $mn$. In the case that $\proj{1} \times \proj{1}$ is
replaced by $\proj{2}$, a similar extension has already been done
by Bus\'e, Cox, and D'Andrea \cite{BCD}. The strategy is to
replace certain of the moving quadrics in the matrix $M$ with the
$k$ linearly independent moving planes which exist because of the
presence of the base points.  The proofs of these results require
an extension of the concept of regularity of a module, which is
traditionally a concept for graded modules, to cover the case of
bigraded modules.  This extension was developed in a recent series
of papers \cite{HW}, \cite{HW2}.  We will start by summarizing the
results needed from these papers, and prove some additional
results needed for the application to our implicitization problem.

\section{Bigraded regularity and Saturation} \label{S:regularity}

We will start by recalling the definition and some of the results
concerning bigraded regularity and saturation as developed in
\cite{HW}.  Our main goal in this section is a bound on the
bigraded regularity of the saturation of a power of an ideal. This
result is inspired by results of Chandler \cite{KC}.

In this section we will work over the polynomial ring $R=K[x_0,\,
\ldots,\, x_m,\, y_0, \ldots,\, y_n]$ where $K$ is an infinite
field, and $m,\,n\ge1$. We will make  $R$ into a bigraded
$K$-algebra in the normal manner by assigning the bidegree
$(1,\,0)$ to the $x_i$ variables   and the bidegree $(0,\,1)$ to
the $y_j$ variables. Moreover, we will partially order $\Z^2$ by
the rule $(k,\,l)\le (r,\,s)$ if $k\le r$ and $l\le s$. As usual,
we let $R_{k,\,l}$ denote the $K$-subspace of $R$ consisting of
bihomogeneous polynomials of bidegree $(k,\,l)$, and if $M$ is a
bigraded $R$-module, then $M_{k,\,l}$ is the $(k,\,l)$
bihomogeneous part of $M$.  Let
$$\mbf{m}=\langle x_0,\,\ldots,\, x_m\rangle\cap
\langle y_0,\, \ldots,\,y_n\rangle =\langle \{x_iy_j: 0\le i\le m,
0\le j\le m \}\rangle \subset R.$$ The ideal $\mbf{m}$ is known as
the \emph{irrelevant ideal} of $R$.  Moreover we note that the
local cohomology modules $H^i_{\mbf{m}}(M)$ are naturally bigraded
$R$-modules.

\begin{Def}\label{D:reg}(See \cite[Definition 3.1]{HW})
We say that a bigraded $R$-module $M$  is \textit{
$(p,p')$-regular} if for all $i\geq 0$,
\[H^i_{\mbf{m}}(M)_{k,k'}=0 \text{ whenever }
(k,k')\in \Reg_{i-1}(p,p'),\] where $\Reg_j(p,p')= \{(x,y)\in
\Z^2: x \geq p-j,\; y \geq p'-j,\; x+y \geq p+p'-j-1\}.$
\end{Def}

\begin{rem} \begin{enumerate}
\item The definition of $(p,\,p')$-regular given here coincides
with  what is called \emph{weakly} $(p,\,p')$-regular in
\cite{HW}. In that paper, a concept known as \emph{strongly}
$(p,\,p')$-regular is also introduced and studied. Since we will
not need this stronger version in this paper, we shall simply use
the term $(p,\,p')$-regular for what would normally be referred to
as weakly $(p,\,p')$-regular.

\item
The definition given above for $(p,p')$-regularity is an extension
to bigraded modules of the concept of Castelnuovo regularity for
graded modules as found, for example, in Ooishi \cite{O}.  The
concept of regularity was originally defined for coherent sheaves
of modules on projective space by Mumford \cite{DM}, and this
version of regularity is also treated in the bigraded case in
\cite{HW}. The reader is referred to this paper for a precise
comparison of the two concepts.  However, in case the bigraded $R$
module $M$ is an ideal $I\subset R$ generated by bihomogeneous
polynomials and $\mathcal{I}\subset \mathcal{O}_X$ (where
$X=\proj{m}\times\proj{n}$) is the corresponding sheaf of ideals
in the structure sheaf $\mathcal{O}_X$, the equivalence is
expressed by the following result. See \cite[Proposition 3.5]{HW}
for details.
\end{enumerate}
\end{rem}
\begin{pro}\label{P:IdealSheafReg} With the above notation,
the bihomogeneous ideal
$I\subset R$ is $(p,\,p')$-regular if and only if natural map
$$I_{p,\,p'}\to H^0(X, \mathcal{I}(p,\,p'))$$
is an isomorphism and $$H^i(X, \, \mathcal{I}(k,\,k')=0 \text{
whenever }  (k,\,k')\in \Reg_i(p,p').$$ As usual,
$\mathcal{I}(k,\,k')$ denotes the twisting of $\mathcal{I}$ in
bidegree $(k,\,k')$. Moreover, if $I$ is $(p,\,p')$-regular, then
the natural map $$I_{d,\,d'}\to H^0(X, \mathcal{I}(d,\,d'))$$ is
an isomorphism for all $(d,\,d')\ge (p,\,p')$.\end{pro}

\begin{Def}\label{D:sat}
Let $M$ be a bigraded submodule of a finitely generated free
$R$-module $F$. The \textit{saturation} of the module $M$, denoted
by $M^{\sat}$ or $\sat(M)$ is the submodule of $F$ defined by
$$M^{\sat}=\{f\in F : \mbf{m}^k f \subset M, \text{ for some }
k\}.$$  The submodule $M$ is said to be \emph{saturated} if
$M=M^{\sat}$, while $M$  is  \textit{$(p,\, p')$-saturated} if
\[ M^{\sat}_{k,\,k'} = M_{k,\,k'}
\text{ for all $(k,\,k')\ge (p,\,p')$.}\]
\end{Def}
\begin{lem}\label{L:H01} Let $R=K[x_0,\,
\ldots,\, x_m,\, y_0, \ldots,\, y_n]$ where $(m,\,n)\ge (1,\,1)$
and let $M$ be a bigraded submodule of a free $R$-module $F$ of
finite rank. Then
\begin{enumerate}
\item $H^0_{\mbf{m}}(M)=0$, and
\item $H^1_{\mbf{m}}(M)\cong M^{\sat}/M$.
\end{enumerate}

\end{lem}
\begin{proof}   $H^0_{\mbf{m}}(M)=
\cup_n(0:_M\mbf{m}^n)=0$ since $M$ is a submodule of a free module
$F$.  The long exact cohomology sequence for
$$\begin{CD} 0 @>>> M @>>> F @>>> F/M @>>> 0\end{CD}$$
has a segment
$$\begin{CD} H^0_{\mbf{m}}(F) @>>> H^0_{\mbf{m}}(F/M)@>>>
H^1_{\mbf{m}}(M) @>>>H^1_{\mbf{m}}(F).\end{CD}$$  Since $F$ is
free and $(m,\,n)\ge(1,\,1)$, it follows that
grade$_F(\mbf{m})\ge2$, so that $H^i_{\mbf{m}}(F)=0$ for $i=0,\,1$
(see \cite[Theorem 6.2.7, Page 109]{BRS}). Thus there is an
isomorphism
$$H^1_{\mbf{m}}(M)\cong H^0_{\mbf{m}}(F/M) =M^{\sat}/M.$$
\end{proof}

\begin{pro}
Let $M$ be a bigraded submodule of a free $R$-module $F$ of finite
rank and let $\mc{M}$ be the corresponding coherent sheaf of
modules on $X=\proj{m}\times\proj{n}$.  Then
$$M^{\sat}_{k,\,l}= H^0(X, \, \mc{M}(k,\,l)).$$
\end{pro}
\begin{proof}For any finitely generated bigraded $R$-module $M$
there is an exact sequence (see \cite[Corollary 1.5]{HY}):
\[\begin{CD} 0 @>>> H^0_{\mbf{m}}(M) @>>> M @>>>
\bigoplus^{}_{(a,\,b)\in \Z^2}H^0(X,\, \mc{M}(a,\,b)) @>>>
H^1_{\mbf{m}}(M) @>>> 0.
\end{CD}
\]
Since $M$ and $M^{\sat}$ generate the same sheaf $\mc{M}$ on $X$,
we can apply this exact sequence with $M$ replaced by $M^{\sat}$.
Lemma \ref{L:H01} shows that $H^i_{\mbf{m}}(M^{\sat})=0$ for
$i=0,\,1$, and the result follows.
\end{proof}
\begin{cor} \label{C:H01}
If $M$ is a bigraded submodule of a free $R$-module $F$ of finite
rank, then $M$ is $(p,\, p')$-saturated if and only if
$H^1_{\mbf{m}}(M)_{k,\,k'}=0$ for all $(k,\,k')\ge(p,\,p')$.
Moreover, if $M$ is $(p,\,p')$-regular, then $M$ is
$(p,\,p')$-saturated.
\end{cor}

The converse of the last statement is true in the case of
dimension 0:
\begin{lem}\label{L:satregzero}
Let $I\subset R$ be a bihomogeneous ideal with $\dim R/I=0$, where
$\dim$ refers to Krull dimension.  Then the following are
equivalent:
\begin{enumerate}
\item $I$ is $(p,\, p')$-saturated.
\item $I$ is $(p,\, p')$-regular.
\item $I_{k,\,k'} = R_{k,\,k'}$ for all $(k,\,k')\ge (p,\,p')$.
\end{enumerate}
\end{lem}

\begin{proof}
$(1. \Leftrightarrow 3.)$ This is clear since $I$ is $\langle
x,y\rangle$-primary.

$(2. \Rightarrow 1.)$ Corollary \ref{C:H01}.

$(1. \Rightarrow 2.)$  According to Definition \ref{D:reg}, we
need to show that $$H^i_{\mbf{m}}(I)_{k,\,k'}=0 \text{ whenever
$(k,\,k')\in \Reg_{i-1}(p,\,p')$.} \leqno{(*)}$$ Since
$H^0_{\mbf{m}}(I)=0$, $(*)$ is certainly true for $i=0$, and since
$I$ is $(p,\,p')$-saturated,
$H^1_{\mbf{m}}(I)_{k,k'}=I^{\sat}_{k,k'}/I_{k,k'}=0$ for all
$(k,k')\in(p,\,p')+\Z_+^2= \Reg_0(p,p')$, where
$\Z_+=\{x\in\Z:x\ge 0\}$. Thus $(*)$ is satisfied for $i=1$. Now
consider the case $i\ge 2$.  The long exact cohomology sequence of
the exact sequence
\[ \begin{CD} 0 @>>> I @>>> R @>>> R/I @>>> 0\end{CD}\]
contains the segment
\[\begin{CD} H^{i-1}_{\mbf{m}} (R/I) @>>> H^i_{\mbf{m}}(I) @>>>
H^i_{\mbf{m}}(R) @>>> H^i_{\mbf{m}}(R/I) \end{CD}.\] Since $\dim
R/I=0$, it follows that $H^i_{\mbf{m}}(R/I)=0$ for $i \geq 1$, so
that if $i \geq 2$, we conclude that
$H^i_{\mbf{m}}(I)=H^i_{\mbf{m}}(R)$.  By \cite[Proposition 4.3 and
Corollary 4.5]{HW}, $R$ is $(0,0)$-regular, and by \cite[Theorem
3.4]{HW}, it follows that $R$ is $(p,\,p')$-regular for all
$(p,\,p')\ge (0,\,0)$. Therefore,
$H^i_{\mbf{m}}(I)_{k,k'}=H^i_{\mbf{m}}(R)_{k,k'}= 0$ for all
$(k,k')\in \Reg_{i-1}(p,\,p')$.  Thus, $(*)$ is also satisfied for
$i\ge 2$, and hence $I$ is $(p,\ p')$-regular.
\end{proof}

With this background out of the way we can proceed with a
discussion of the results on regularity of the powers of a
bihomogeneous ideal that will be needed for the implicitization
problem.

\begin{pro}
\label{P:powerzero} Let $I \subset R$ generated by bihomogeneous
forms of bidegree $ \le (r,\,r')$, and assume that $I$ is
$(p,\,p')$-regular. If $\dim R/I=0$, then $I^e$ is
$(l,\,l')$-regular for some $(l,\,l') \leq ((e-1)r+p,\,
(e-1)r'+p')$.
\end{pro}

\begin{proof}
The proof is by induction on $e$. The result is true for $e=1$ by
assumption. Since $\dim R/I=\dim R/I^e=0$,  we can proceed by
induction, and assume that $I^{e-1}$ is
$((e-2)r+p,\,(e-2)r'+p')$-regular.

According to Lemma \ref{L:satregzero}, we need to show that
$I^e_{k,k'}=R_{k,k'}$ for any $(k,\,k')\geq ((e-1)r+p,\,
(e-1)r'+p')$.  For this, it will suffice to show that $M \in I^e$
for every monomial $M$ of bidegree $(k,\,k')$, where $(k,\,k')\geq
((e-1)r+p,\, (e-1)r'+p')$.  Thus let $M$ be an arbitrary monomial
of bidegree $(k,\,k')\geq ((e-1)r+p,\, (e-1)r'+p')$. Write $M$ as
a product $M=NN'$ where $N$ and $N'$ are monomials of bidegrees
$(p,\,p')$ and $(k-p,\,k'-p')$, respectively. Suppose that
$I=\langle f_1, \ldots, f_s\rangle$, where $\bideg(f_i)=(d_i, \,
d'_i)\le (r,\,r')$, for all $i$. Since $I$ is $(p,\, p')$-regular,
$N\in I$ by Lemma \ref{L:satregzero}. Thus we can write
$N=\sum_{i=1}^{s}N_if_i$, where $\bideg(N_i)=(n_i,
n'_i)=(p-d_i,\,p'-d_i')\ge(p-r,\,p'-r')$, so that
$$\bideg(N_iN')\ge (k-r,\,k-r')\ge ((e-2)r+p,\,(e-2)r'+p').$$
By the induction hypothesis, $N_iN' \in I^{e-1}$, and hence
$M=\sum_{i=1}^r N_iN'f_i \in I^{e-1}I=I^e$. By Lemma
\ref{L:satregzero}, we conclude that $I^e$ is $((e-1)m+p,\,
(e-1)m'+p')$-regular.
\end{proof}

\begin{thm}\label{T:satpower}
Let $I \subset R$ be a bihomogeneous ideal, and  assume that
\begin{itemize}
\item[1.] $\V(I)\subset \proj{m}\times \proj{n}$ is finite;
\item[2.] $I$ is  $(p,p')$-regular;
\item[3.] $I$ is generated by forms of bidegree $\le (r, \, r')$.
\end{itemize}
Then $J=\sat(I^e)$ is  $((e-1)r+p, \,(e-1)r'+p')$-regular.
\end{thm}

\begin{proof} We will let $X=\proj{m}\times \proj{n}$ and $Z$ will
denote the finite subscheme $\V(I)$. The proof is by induction on
$e$. Suppose $e=1$. In this case, it is necessary to show $J$ is
$(p,\, p')$-regular, i.e., $H^i_{\mbf{m}}(J)_{k,k'}=0$ for all
$i\ge 0$ and  $(k,\, k')\in \Reg_{i-1}(p,\,p')$. Since $J$ is a
saturated ideal, Lemma \ref{L:H01} shows that $H^i_{\mbf{m}}(J)=0$
for $i=0, 1$, so that $H^i_{\mbf{m}}(J)_{k,k'}=0$ for $i=0,1$ and
for all $k,\,k'$.

If $i\geq 2$, let $\mc{I}$ and $ \mc{J}$ be the sheaves on
$X=\proj{m}\times \proj{n}$ defined by $I$, and $J$ respectively.
The long exact cohomology sequence of the  exact sequence
\[ \begin{CD} 0 @>>> \mc{I} @>>> \mc{J} @>>> \mc{J}/\mc{I}
@>>> 0 \end{CD} \] tensored with $\mc{O}(k,\,k')$ contains the
segment
\[ H^{i-1}(Z, (\mc{J}/\mc{I})(k,\,k')) \rightarrow H^i(X,
\mc{I}(k,\,k')) \rightarrow H^i(X, \mc{J}(k,\,k')) \rightarrow
H^i(Z, (\mc{J}/\mc{I})(k,\,k')) .\] Since $\dim Z =0$, $H^i(Z,
(\mc{J}/\mc{I})(k,k'))=0$ for $i\geq 1$ and for all $k,\,k'$.
Since $I$ is $(p,\, p')$-regular, $H^1(X,\,\mc{I}(k,\,k'))=0$ for
$(k,\, k')\in \Reg_1(p,\, p')$ by Proposition
\ref{P:IdealSheafReg}. Thus, we have
\begin{equation}\label{E:only1}
H^1(X, \mc{J}(k,k'))=0, \   \ \forall (k,\,k')\in \Reg_1(p,\,p'),
\end{equation}
and,
\[H^i(X, \mc{J}(k,k'))=H^i(X, \mc{I}(k,k')), \, \ \forall i\geq 2 .\]
Since $I$ is $(p,\,p')$-regular,
\begin{equation}\label{E:atleast2}
H^i(X, \mc{J}(k,k'))=H^i(X, \mc{I}(k,k'))=0, \  \ \forall i\geq 2
,\  \ \forall (k,\,k')\in \Reg_i(p,\,p'),
\end{equation}
and combining Equations \eqref{E:atleast2} and \eqref{E:only1}, we
conclude that
\[H^i(X, \mc{J}(k,k'))=0, \ \ \forall
i\geq 1, \  \ (k,\,k')\in \Reg_i(p,\,p').\] Since
$H^{i+1}_{\mbf{m}}(J)_{k,k'}=H^i(\mc{J}(k,k'))$,  for all $i \geq
1$, it follows that
\[H^{i}_{\mbf{m}}(J)_{k,k'}=0, \  \ \forall
i\geq 2, \  \ (k,k')\in \Reg_{i-1}(p,p'),\] and hence, $J$ is
$(p,p')$-regular when $e=1$.

Now assume that $e \ge 2$. The sheafification of $J$ is
$\mc{I}^e$, and $H^0(X, \mc{I}^e(k,k'))=J_{k,k'}$. Define
$Z^{(d)}=\V(I^d )$, which has the same support as $Z$ and is hence
is finite. Since $J$ is saturated, we have $H^i_{\mbf{m}}(J)=0$
for $i=0,1$ (Lemma \ref{L:H01}). Let $(l, l')
=((e-1)r+p,\,(e-1)r'+p')) $. We must show that
\[
H^i(X, \mc{I}^e(k,k'))=0 \text{ for } (k,k')\in \Reg_i(l,l'),
\text{ all } i \ge 1.
\]
Tensor the following exact sequence
\[
\begin{CD} 0@>>>\mc{I}^e @>>> \mc{I}^{e-1} @>>>
\mc{I}^{e-1}/ \mc{I}^e @>>>0.
\end{CD}
\]
with $\mc{O}(k,k')$ and consider the resulting cohomology
sequence. Since the support of $\mc{I}^{e-1}/ \mc{I}^e$ is
contained in $Z$, which is $0$-dimensional,  it follows that
$\coh{i}{X}{(\mc{I}^{e-1} /\mc{I}^{e})(k, k')} =0$ for $i\geq 1$.
Therefore,
\[
\coh{i}{X}{\mc{I}^{e} (k,\, k')} = \coh{i}{X}{\mc{I}^{e-1} (k,
k')} \text{ for all $i \ge 2$},
 \]
and the latter group vanishes by induction for all
\[
(k,k')\in \Reg_i((e-2)r+p, \,(e-2)r'+p')\supset
 \Reg_i((e-1)r+p,\,(e-1)r'+p').
\]
Thus, we have the required vanishing for $i \ge 2$. Now look at
the sequence
\[{\small
\begin{CD}
H^0(X, \mc{I}^{e-1}(k,k'))@>\phi>> H^0(X, (\mc{I}^{e-1}/
\mc{I}^e)(k,k')) @>>>H^1(X, \mc{I}^{e}(k,k')) @>>>H^1(X,
\mc{I}^{e-1}(k,k'))
\end{CD}}
\]
By induction, the last term vanishes for all $(k, k')\in \Reg_1
(l, l')$, so that the next-to-last term will vanish there provided
we show that $\phi $ is onto for those same $(k, k')$.
\par
Suppose $Z=\{\mbf{p}_1, \ldots, \mbf{p}_s\}$. Note, since the
support is finite, we have
\[H^0(X, (\mc{I}^{e-1}/ \mc{I}^e)(k,k'))=
\bigoplus_{\mbf{p}\in Z } (I^{e-1}\mc{O}_{X,\mbf{p}}/I^e
\mc{O}_{X,\mbf{p}})(k,k')
\]
We will show that for $(k, k')\in \Reg_1 (l, l')$ and for any
\[
\left ( \frac{u _1}{v _1}, \dots,  \frac{u _s}{v _s}\right ) \in
\bigoplus_{i} (I^{e-1}\mc{O}_{X,\mbf{p}_i}/I^{e}
\mc{O}_{X,\mbf{p}_i})(k,k')
\]
with bihomogeneous forms with $\deg u_i - \deg v_i = (k, k')$,
$u_i \in I^{e-1}$ we can find a bihomogeneous $g\in
(I^{e-1})^{\sat}_{k,k'}$ and forms $H_i$ with $H_i (\mbf{p}_i )
\ne 0$, such that
\begin{equation}
\label{E:power0} H_i (g v_i - u_i) \in I^{e} \text { for all } i.
\end{equation}
This will prove that  $\phi$ is surjective.

Let $I$ be generated by bihomogeneous elements $f_1, \dots, f_r$
with bidegree $(m_i, \, m_i')\le (r,\,  r')$. We can write
\[
u _i = \sum a_{ij}f_j, \text {for some } a_{ij}\in I^{e-2}_{k-m_j,
\, k'-m_j'}
\]
Note that $(\alpha,\, \alpha') = (k-m_j,\, k'-m_j') \in \Reg_1
((e-2)r +p,\, (e-2)r' +p')$, by our initial choice of $(k,\,k')$.
Tensor the following exact sequence
\[
\begin{CD}
0 @>>> \mc{I}^{e-1} @>>> \mc{I}^{e-2} @>>> \mc{I}^{e-2}/
\mc{I}^{e-1} @>>>0
\end{CD}
\]
with $\mc{O}_X(\alpha,\alpha')$, and consider the resulting
cohomology sequence
\[
\begin{CD}
H^0(X, \mc{I}^{e-2}(\alpha,\alpha'))@>\psi>> H^0(X, (\mc{I}^{e-2}
/\mc{I}^{e-1}) (\alpha, \alpha'))\end{CD}\]
\[\begin{CD}
@>>>H^1(X, \mc{I}^{e-1}(\alpha,\alpha')) @>>>H^1(X,
\mc{I}^{e-2}(\alpha,\alpha')).
\end{CD}
\]
By our induction hypothesis, the third term vanishes, so that
$\psi$ is onto for this $(\alpha,\alpha')$. This means that for
every $j$, and each
\[
\left ( \frac{a _{1j}}{v _1}, \dots, \frac{a _{sj}}{v _s}\right )
\in \bigoplus_{i}( I^{e-2}\mc{O}_{X,\mbf{p}_i}/I^{e-1}
\mc{O}_{X,\mbf{p}i})(k-m_j,k'-m'_j)
\]
we can find a bihomogeneous $g_j\in
(I^{e-2})^{\sat}_{\alpha,\alpha'}$ and forms $H_{ij}$ with $H_{ij}
(\mbf{p}_i ) \ne 0$, such that
\begin{equation}
\label{E:power} H_{ij} (g_j v_{i} - a_{ij}) \in I^{e-1} \text {
for all } i.
\end{equation}
We may replace each $H_{ij}$ by $ H_i = \prod _j H_{ij}$. Multiply
equation (\ref{E:power}) by $f_j$ and sum the result over $j$ and
define $g = \sum g_j f_j \in (I^{e-1})^{\sat}_{k,k'}$. Then we
have obtained equation (\ref{E:power0}), as required.
\end{proof}

\section{Finite Subschemes of $\proj{1}\times\proj{1}$}\label{S:finite}

This section will be devoted to a presentation of several results
which can be proven for bihomogeneous ideals $I$ that define
finite subschemes of $\proj{1}\times\proj{1}$, but that do not
necessarily have immediate analogues for subschemes of general
biprojective spaces. The results proved are analogous to the
results proved in \cite{BCD} for application to the
implicitization problem for maps $\phi:\proj{2}\to \proj{3}$. Our
results will be similarly applied for the implicitization of maps
$\phi:\proj{1}\times \proj{1}\to \proj{3}$.  Since we are
restricting ourselves to this low dimensional case, we will let
$R$ be the polynomial ring $\C[s,\,u,\,t,\,v]$ in the variables
$s$, $u$, $t$, and $v$, where, as usual, the bigrading of $R$ is
given by setting the bidegree of $s$ and $u$ to be $(1,\,0)$ and
the bidegree of $t$ and $v$ to be $(0,\,1)$. If  $I=\langle f_1,
\dots, f_r \rangle\subset R $ is an ideal generated by forms all
of the same bidegree $(m,n)$ with $m, n \geq 1$, then there is a
rational map $\phi_I:\proj{1}\times\proj{1}\to \proj{r-1}$ defined
by
$$ \phi_I =[f_1(s,\,u;t,\,v),\ldots,f_r(s,\,u;t,\,v)].$$

Note that the polynomial ring $\C[s,\,t,\,v]$ inherits a bigrading
as a subring of $R=\C[s,\,u,\,t,\,v]$, so that a polynomial
$f(s,\,t,\,v)$ is bihomogeneous with bidegree $(m,\,n)$ if and
only if $f(s,\,t,\,v)=\sum_{j=0}^n
a_{ij}s^mt^iv^{n-j}=s^mg(t,\,v)$, where $g(t,\,v)$ is homogeneous
of degree $n$.
\begin{lem}\label{L:IJ}
Let $\bar{I} \subset S=\C[s,\,t,\,v]$ be an ideal,  minimally
generated by $r$ bihomogeneous forms of bidegree $(m,\,n)$.  That
is, $\bar{I}=s^m J$ where $J\subset \C[t,\,v]$ is generated by
homogeneous polynomials of degree $n$. If $\V(J)=\emptyset$ in
$\proj{1}$, then $\bar{I}$ is $(p,\,p')$-regular for all $p\geq m$
and $p' \geq 2n-r+1$.
\end{lem}
\begin{proof}
This follows from \cite[Remark 4.12]{HW} and Lemma B.1 in
\cite{BCD}.
\end{proof}

\begin{rem}
Similarly, let $\bar{I} \subset S=\C[s,u,t]$ be an ideal,
minimally generated by $r$ bihomogeneous forms of bidegree
$(m,n)$. That is $\bar{I}=t^n J$ where $J$ is generated by
homogeneous polynomials in $\C[s,\,u]$ of degree $m$. If
$\V(J)=\emptyset$ in $\proj{1}$, then $\bar{I}$ is
$(p,p')$-regular for all $p \geq 2m-r+1$ and $p'\geq n$.
\end{rem}

\begin{lem}\label{L:mingen}
Let $I \subset R=\C[s,\,u,\,t,\,v]$ be minimally  generated by $r
\ge 4$ bihomogeneous forms of bidegree $(m,\,n)$ with both $m,\ n
\geq 1$. Assume that $\dim \im{\phi_I}=2$ and that $\V(I) \subset
\proj{1}\times \proj{1}$ is finite. Given $\ell \in R_{1,\,0}$,
let $I_\ell$ be the image of $I$ in the quotient ring $R/\langle
\ell \rangle$.  Then for a generic $\ell$,  $I_\ell $ is minimally
generated by at least 2 elements.
\end{lem}

\begin{proof} The proof is a straightforward modification of the
Bertini theorem argument in \cite[Lemma B.2]{BCD}. See also
\cite[Lemma 3.4.3]{Wang}.
\end{proof}

\begin{rem}
The above result is also true if the given generic element $\ell$
is chosen from $R_{0,\,1}$.
\end{rem}

The following is the main vanishing theorem needed for our
applications.
\begin{thm}\label{T:degbd}
Let $I\subset R=\C[s,u,t,v]$ be minimally generated by $r \ge 4$
bihomogeneous forms of bidegree $(m,\,n)$. Assume that $\dim
\im{\phi_I}=2$ and  assume that $\V(I) \subset \proj{1} \times
\proj{1}$ is finite. If $\mc{I}$ is the associated sheaf of ideals
on $X= \proj{1}\times \proj{1}$, then
\begin{enumerate}
\item $H^1(X,\,\mc{I}(k,k'))=0$ for all $(k,\,k') \geq (2m-2,\,2n-2)$,
and
\item $H^2(X,\,\mc{I}(k,k'))=0$ for all $(k,\, k') \geq (0,\,0)$.
\end{enumerate}
\end{thm}

\begin{proof}
If $Z=\V(I) \subset X=\proj{1} \times \proj{1}$, there is an exact
sequence
\[ 0 \rightarrow \mc{I} \rightarrow
\mc{O}_{X} \rightarrow \mc{O}_Z \rightarrow 0,\] which, upon
taking the tensor product with $\mc{O}_X(k,k')$, gives rise to a
long exact cohomology sequence
\[ \rightarrow H^1(X,\,\mc{O}_Z(k,k')) \rightarrow H^2(X,\,\mc{I}(k,k'))
\rightarrow H^2(X,\,\mc{O}_X(k,k')) \rightarrow
H^2(X,\,\mc{O}_Z(k,k')) \rightarrow .\] Since $Z$ is finite,
$H^i(X,\,\mc{O}_Z(k,\,k'))=0$ for all $i \geq 1$.  By the
K\"{u}nneth formula \cite{SW},
\[ H^2(X,\,\mc{I}(k,k'))= H^2(X,\,\mc{O}_{X}(k,k'))=0,
 \text{ for all  $(k,\,k') \geq (0,\,0)$}. \]
This proves item 2.

To prove the first statement, choose a line $\ell \in R_{1,\,0}$
such that $\V(\ell) \cap \V(I) = \emptyset $ and
$\bar{I}=I_{\ell}=$ the image of $I$ in $R/ \langle \ell \rangle$
is minimally generated by at least two elements. This is possible
by Lemma \ref{L:mingen}. Then by Lemma \ref{L:IJ}, we know that
$\bar{I}$ is $(p,\,p')$-regular for $p \geq m$ and $p'\geq 2n-1$.
If $\bar{\mc{I}}$ is the sheaf on $\V(\ell)\cong \proj{1}$
associated to $\bar{I}$, then by Proposition
\ref{P:IdealSheafReg}, we have
\begin{equation}\label{E:H01}
\begin{array}{rl}
\bar{I}_{k,k'} \cong H^0(\V(\ell),\,\bar{\mc{I}}(k,k')) &\text{
for
all $(k,\,k') \geq (m,\,2n-1)$, and}\\
H^1(\V(\ell),\,\bar{\mc{I}}(k,k'))=0& \text{ for all $(k, \, k')
\geq (m-1,
 \, 2n-2)$.}
 \end{array}
 \end{equation}
Now, we consider the following exact sequence:
\[0  \rightarrow \mc{O}_{\proj{1}\times\proj{1}}(-1,0) \rightarrow
\mc{O}_{\proj{1}\times \proj{1}} \rightarrow \mc{O}_{\V(\ell)}
\cong \mc{O}_{\proj{1}}
\rightarrow 0.\] Tensoring with $\mc{I}(k,k')$ gives the exact
sequence:
\[ Tor_1^{\mc{O}_{\proj{1}\times\proj{1}}}(\mc{I}(k,k'),
\mc{O}_{\proj{1}}) \rightarrow \mc{I}(k-1,k') \rightarrow
\mc{I}(k,k') \rightarrow \mc{O}_{\proj{1}}
\otimes_{\mc{O}_{\proj{1}\times \proj{1}}} \mc{I}(k,k')
\rightarrow 0. \]
 Note $\mc{O}_{\proj{1}}
\otimes_{\mc{O}_{\proj{1}\times \proj{1}}} \mc{I}(k,k') \cong
\bar{\mc{I}}(k,k')$. Since $\mc{O}_{\proj{1}}$ is supported on
$\V(\ell)$, the $Tor$-sheaf is supported there.  Also, for
$p\notin \V(I)$, the sheaf $\mc{I}(k,k')$ is locally free.  Hence
the $Tor$-sheaf  vanishes at $p$ if $p\notin\V(I)$. Hence the
support of the $Tor$-sheaf is contained in $\V(I) \cap V(\ell)$.
By the generic choice of $\ell$, $\V(I)\cap \V(\ell)=\emptyset$,
so the $Tor$-sheaf vanishes. Thus there is exact sheaf sequence
\[0 \rightarrow \mc{I}(k-1,k') \rightarrow
\mc{I}(k,k') \rightarrow  \bar{\mc{I}}(k,k') \rightarrow 0. \]
that gives the following commutative diagram
\[\begin{matrix}\begin{CD}
  &I_{k,k'} & @>>>& \bar{I}_{k,k'} &@>>> & 0  \\
 & \downarrow &  & & &     \downarrow    & &  &   &\downarrow  \\
 & H^0(X,\,\mc{I}(k,k')) & @> \alpha>> &
H^0(\V(\ell),\,\bar{\mc{I}}(k,k')) & @> \beta>> &
H^1(X,\,\mc{I}(k-1,k'))
\end{CD} \end{matrix}\]
\[\begin{CD} @>>> H^1(X,\,\mc{I}(k,k')) @>>> H^1(\V(\ell),\,
\bar{\mc{I}}(k,k'))
\end{CD} \] with exact rows.
If $(k,\,k')\geq (m,\,2n-1)$, Equation \eqref{E:H01} shows that
$H^1(\V(\ell),\,\bar{\mc{I}}(k,k'))=0$ and $\bar{I}_{k,k'} \cong
H^0(\V(\ell),\,\bar{\mc{I}}(k,k'))$. Therefore,  $\alpha$ is onto,
and  $\beta$ is  zero, which implies that there is an isomorphism
\[ H^1(X,\,\mc{I}(k-1,k')) \cong H^1(X,\,\mc{I}(k,k')), \text{  for all
$ (k,\,k')\ge (m,\, 2n-1)$}.\]

An analogous argument with a generic line $\ell \in R_{0,\,1}$
produces another isomorphism
\[ H^1(X,\,\mc{I}(k,k'-1)) \cong H^1(X,\,\mc{I}(k,k')), \text{  for all
$ (k,\,k')\ge (2m-1,\, n)$}.\] Therefore,
\[H^1(X,\,\mc{I}(k-1,k'-1)) \cong H^1(X,\,\mc{I}(k,k')),\text{  for all
$ (k,\,k')\ge (2m-1,\, 2n-1)$}.\] Since $H^1(X,\,\mc{I}(m,n))=0$
if $(m, \, n) \gg (0,\,0)$, we conclude that
$H^1(X,\,\mc{I}(k,k'))=0$  for all $( k ,\, k')\ge (2m-2,\,
 2n-2)$.
\end{proof}

We are now able to prove the following result relating regularity
of the ideal $I$ and the degree of the 0-dimensional subscheme
$\V(I)$. This is one of the main results needed for the
application to the implicitization problem.

\begin{thm}\label{T:wkreg} Let $I\subset R$ be minimally generated by
$r\geq 4$ bihomogeneous forms of bidegree $(m,n)$ with $m,n\geq
1$. Assume that $\mathbb{V}(I) \subset X=\proj{1} \times \proj{1}$
is finite and $\dim \im{\phi_I}=2$. If $(p,\,p')\ge
(2m-1,\,2n-1)$, then $I$ is  $(p,p')$-regular if and only if
$\dim_{\C}(R/I)_{p,p'}=\deg(\V(I))$, where $\deg(\V(I))$ denotes
the degree of the $0$-dimensional subscheme $\V(I)$.
\end{thm}

\begin{proof}
 When $p\geq 2m-1$ and $p' \geq
2n-1$, Theorem \ref{T:degbd} implies $H^1(X,\,\mc{I}(p,p'))=0$.
Thus, the exact sheaf sequence
\[ 0 \rightarrow \mc{I} \rightarrow
\mc{O}_{\proj{1} \times \proj{1}} \rightarrow \mc{O}_Z \rightarrow
0\] produces an exact sequence
\[ 0 \rightarrow H^0(X,\,\mc{I}(p,p')) \rightarrow
H^0(X,\, \mc{O}_{X}(p,p')) \rightarrow H^0(Z,\,\mc{O}_Z(p,p'))
\rightarrow 0.\] This gives the following commutative diagram with
exact rows:
\[\begin{matrix}
 0 & \rightarrow & H^0(X,\,\mc{I}(p,p'))& \rightarrow &
 H^0(X,\,\mc{O}_{X}(p,p'))&
\rightarrow & H^0(Z,\,\mc{O}_Z(p,p'))& \rightarrow &  0 \\

  &              & \downarrow   &    & \downarrow  &   &
  \downarrow &  & \\
0 & \rightarrow & I_{p,p'} & \rightarrow & R_{p,p'} & \rightarrow
& (R/I)_{p,p'} & \rightarrow & 0.\end{matrix}\] We have
$R_{p,p'}=H^0(X,\,\mc{O}_{X}(p,p'))$ and if $I$ is
$(p,p')$-regular, then $I_{p,p'}=H^0(X,\,\mc{I}(p,p'))$. The
5-lemma then shows that $(R/I)_{p,p'}=H^0(Z,\,\mc{O}_Z(p,p'))$, so
that $$\dim_{\C}(R/I)_{p,p'}=\dim_{\C} H^0(Z,\,\mc{O}_Z(p,p')).$$
But
\[\dim_{\C}H^0(Z,\,\mc{O}_Z)=\deg(Z),\]
and since $$\dim_{\C}H^0(Z,\,\mc{O}_Z) = \dim_{\C}
H^0(Z,\,\mc{O}_Z(p,p'))$$ when $Z$ is finite, we conclude that
\[\dim_{\C}(R/I)_{p,p'}=\deg(Z).\]

Conversely, suppose $\dim_{\C}(R/I)_{p,p'}=\deg(Z)$. Since
$H^2(X,\,\mc{I}(k,k'))=0$ for all $k,k' \geq 0$ by  Theorem
\ref{T:degbd}, it follows from Proposition \ref{P:IdealSheafReg}
that to show  $I$ is $(p,p')$-regular, we only need to prove that
\[ I_{p,p'} \cong H^0(X,\,\mc{I}(p,p')), \text{ and }
H^1(X,\,\mc{I}(p-1,p'-1))=0 . \] If $p\geq 2m-1$ and $p'\geq
2n-1$, then $H^1(X,\,\mc{I}(p-1,p'-1))=0$ by Theorem
\ref{T:degbd}. We know that the natural map $I_{p,p'} \rightarrow
H^0(X,\,\mc{I}(p,p'))$ is injective, so it is enough to show that
$$\dim_{\C} I_{p,p'}=\dim_{\C} H^0(X,\,\mc{I}(p,p')).$$ From
the exact sequence
\[ 0 \rightarrow H^0(X,\,\mc{I}(p,p')) \rightarrow R_{p,p'}
\rightarrow H^0(Z,\,\mc{O}_Z(p,p')) \rightarrow 0, \] we conclude
that
\[\dim_{\C} H^0(X,\,\mc{I}(p,p'))=\dim_{\C} R_{p,p'} -
\dim_{\C} H^0(Z,\,\mc{O}_Z(p,p'))\]
\[=\dim_{\C} R_{p,p'} -\deg(Z) =\dim_{\C} R_{p,p'} -
\dim_{\C} (R/I)_{p,p'}=\dim_{\C} I_{p,p'}.\] Thus $I_{p,p'} \cong
H^0(X,\,\mc{I}_{p,p'})$ and $I$ is $(p,p')$-regular.
\end{proof}

\begin{cor}\label{L:wkregdim0}
Under the hypotheses of Theorem \ref{T:wkreg}, if $I$ is
$(p,p')$-regular, then $\dim_{\C}(R/I)_{k,k'}=\deg(\V(I))$ for all
$(k,k')\geq (p,p')$.
\end{cor}
\begin{proof}
If $I$ is $(p,\,p')$-regular, then $I$ is $(k,\,k')$-regular for
all $(k,\,k')\ge (p,\,p')$.
\end{proof}

\begin{exa}
If $I=\langle u^2t^2v, u^2t^3+suv^3, s^2tv^2,
s^2v^3+s^2t^3\rangle\subset \C[s,\,u,\,t,\,v]$, then
$\V(I)=(0,1;0,1)\in \proj{1}\times \proj{1}$. In this case, each
generator of $I$ has bidegree $(m,n)=(2,3)$ and
$(2m-1,2n-1)=(3,5)$. A computation with Singular \cite{GPS01}
shows $\dim_{\C}(R/I)_{3,5}=\deg \V(I)=2$. Therefore, $I$ is
$(3,5)$-regular by Theorem \ref{T:wkreg}.
\end{exa}

We will conclude this section with a brief description of a result
on syzygies that will be needed in the proof of our
implicitization theorem.

\begin{Def}
\label{D:base} Let $I=\langle r_1,\,\ldots,\, r_n\rangle \subseteq
R$ be an ideal generated by bihomogeneous elements of $R$. In
analogy with the case of a rational map, we will say that
$\mathbb{V}(I)$ is the {\it base point scheme} of $I$.
\begin{enumerate}
\item The \emph{syzygy submodule of $I$} is the submodule of
relations among the $r_i$ ($1\le i\le n$) defined by
$$ \Syz(r_1,\, \ldots,\, r_n) = \set{(a_1,\, \ldots,\, a_n)\in
R^n: a_1r_1+\cdots+a_nr_n=0}.$$
\item
A syzygy $(a_1, \ldots, a_n)\in \Syz(r_1, \ldots, r_n)$
\textit{vanishes at the base points of} $I$ if, for each $i$,
$a_i\in I^{\sat}$.

\item
A syzygy $(a_1, \ldots, a_n)\in \Syz(r_1, \ldots, r_n)$ has
bidegree $(k,\,l)$ provided each $a_i$ has bidegree $(k,\,l)$.

\item If $\mathbf{e}_i\in R^n$ denotes the standard basis vector
with a $1$ in the $i^{\rm th}$ position and 0 elsewhere, then a
\emph{basic Koszul syzygy} for the ideal $I$ is one of the form
$$\mathbf{s}_{ij}=r_j\mathbf{e}_i-r_i\mathbf{e}_j, \quad
\text{for $i<j$}.$$ Since $(r_j)r_i+(-r_i)r_j=0$ for $i \neq j$,
it is clear that $\mathbf{s}_{ij}\in \Syz(r_1,\, \ldots,\, r_n)$.
Let $\Kos(r_1,\, \ldots,\, r_n)\subset \Syz(r_1,\, \ldots,\, r_n)$
be the submodule generated by the basic Koszul syzygies. We refer
to an arbitrary element of $\Kos(r_1,\, \ldots,\, r_n)$ as a
\emph{Koszul syzygy}.
\end{enumerate}
\end{Def}

The following result is the fundamental result relating the Koszul
syzygies of the ideal $I$  and the module of syzygies of $I$ which
vanish at the base points of $I$ in the special situation that we
will need for this paper.

\begin{thm}\label{T:lcikoz}
Let $a_0$, $a_1$ and $a_2\in R$ be bihomogeneous polynomials of
bidegree $(m,\,n)$ and suppose that $\V(a_0,\,a_1,\,a_2)\subset
\proj{1}\times \proj{1}$ is finite, each base point $\mbf{p}\in
\V(a_0,\,a_1,\,a_2)$ is a local complete intersection, and that
$(A_0,\,A_1,\,A_2)\in \Syz(a_0,\,a_1,\,a_2)$ is a syzygy of
bidegree $(k,\,l)$, where $(k-2m+1)(l-2n+1)\ge 0$. Then
$(A_0,\,A_1,\,A_2)$ vanishes on the base points of $I=\langle
a_0,\, a_1,\, a_2\rangle$, if and only if $(A_0,\,A_1,\,A_2)\in
\Kos(a_0,\,a_1,\,a_2)$, which means that there are $h_1$, $h_2$,
$h_3$ of bidegree $(k-m,\, l-n)$ such that
\begin{eqnarray*}
A_0&=&h_1a_2+h_2a_1\\
A_1&=&-h_2a_0+h_3a_2\\
A_2&=&-h_1a_0-h_3a_1.
\end{eqnarray*}
\end{thm}

\begin{proof}
This result follows from Corollary 3.15 of \cite{HW2}.  See Remark
3.16 in that paper.
\end{proof}

To say that $I$ is a local complete intersection means that each
local ring $\mathcal{I}_{\mathbf{p}}$ of the associated sheaf of
ideals $\mathcal{I}$  is a complete intersection ideal. Precisely,
if $I$ is an ideal of $R$ generated by bihomogeneous forms, and
$Z=\V(I)\subset \proj{1}\times \proj{1}$ is a finite set, then we
say that a base point $\mathbf{p}\in Z$ is a \textit{local
complete intersection} (LCI) if the local ring
$\mathcal{I}_{\mathbf{p}}\subset \mathcal{O}_{X,\mathbf{p}}$ is a
complete intersection ideal, i.e., $\mathcal{I}_{\mathbf{p}}$ is
generated by two elements. The ideal $I$ is a \emph{local complete
intersection} provided each base point $\mathbf{p}\in Z$ is a
local complete intersection.

\section{Local complete intersection base points of total
multiplicity $k\leq mn$} \label{S:multiplebp}

In this section,  we will extend the method of moving quadrics to
the case where multiple base points are present. Throughout this
section, $\phi$ will be a map $ \phi: \proj{1} \times \proj{1} \to
\proj{3}$ given by $\phi(s,u; t,v) =[a_0,\, a_1,\, a_2,\,a_3]$
where each $a_i\in R$ is a bihomogeneous polynomial of bidegree
$(m,n)$, and $I=\langle a_0,\, a_1,\, a_2,\, a_3 \rangle$. For
convenience of reference, we will list some conditions on $\phi$
related to  base points.  Some of these conditions will be needed
in each of the results of this paper, and they will be referred to
by number (B1 -- B6) as needed.

\begin{enumerate}
\item [B1:] The polynomials $a_i(s,u,t,v)$ ($0\le i\le 3$)
are bihomogeneous of bidegree $m, n$ and are linearly independent
over $\C$.

\item  [B2:] The base point scheme $\V(I)$
consists of a finite number of base points with total multiplicity
$k \leq mn$.

\item [B3:] Each  base point $\mathbf{p}\in \V(I)$ is a LCI.

\item [B4:] $\dim_{\C}(R/I)_{2m-1,\,2n-1} =\deg (\V(I))$.

\item [B5:]  The base point scheme $\V(I)=\V(a_0,\,a_1,\,a_2)$
and $a_3 \in \sat\langle a_0, a_1, a_2\rangle$.

\item [B6:] $\dim_{\C}
\Syz(a_0,\,a_1,\,a_2)_{m-1,\,n-1}=0$.
\end{enumerate}

\begin{rem}
Some remarks concerning these conditions:
\begin{enumerate}
\item The condition B1 simply says that $S=\im \phi$ is not contained
in any plane in $\proj{3}$.
\item
The finiteness of $\V(I)$ in condition B2 is equivalent to
$\gcd(a_0, a_1, a_2, a_3)=1$, while $k \leq mn$ is equivalent to
the degree inequality $ \deg S \deg \phi \geq mn$. The last
equivalence is a consequence of  the degree formula
\[2mn=\deg \phi \deg S +\sum_{p\in \V(I)}e(I,\mathbf{p}),\]
which is similar to \cite[Page 19]{C}.  For a proof, see
\cite[Theorem 4.2.12]{Wang}.  In this formula, $e(I,\mathbf{p})$
is the multiplicity of the local ring $\mathcal{I}_{\mathbf{p}}$.

\item
The above degree formula for the image of the parametrization
involves the sum of the multiplicities of the base points.  This
equals $\deg(\V(I))$ only when $\V(I)$ is a local complete
intersection.  Hence the need for the condition B3.
\item
Condition B4 is necessary to be able to apply the regularity
condition on $I$ given by Theorem \ref{T:wkreg}.
\item
Conditions B5 and B6 are technical conditions which are needed to
be able to apply Theorem \ref{T:lcikoz}.
\end{enumerate}
\end{rem}

\begin{lem}\label{R:satabc}
Suppose $a_0,\, a_1,\, a_2,\, a_3 \in \C[s,\,u,\,t,\,v]$ are
bihomogeneous of bidegree $m,\,n$ with no common factor, and
$\V(a_0,\, a_1,\, a_2,\, a_3)$ is a local complete intersection.
If we replace $\{a_i\}_{i=0}^2$ with generic linear combinations
of $\{a_i\}_{i=0}^3$, then we have $\V(a_0,\, a_1,\, a_2)=
\V(a_0,\, a_1,\, a_2,\, a_3)$ as subschemes of $\proj{1}\times
\proj{1}$, and $a_3\in \sat\langle a_0,\,a_1,\,a_2\rangle$.
\end{lem}
\begin{proof}
The result is proved in \cite[Theorem A.1, Corollary A.2]{BCD} for
the case  of homogeneous polynomials in $k[x,\,y,\,z]$, but the
argument works verbatim in the case of bihomogeneous polynomials.
\end{proof}

\begin{rem}
A consequence of Lemma \ref{R:satabc} is that if the
parametrization $\phi:\proj{1}\times \proj{1}\to \proj{3}$
satisfies conditions B1 -- B4, then after a generic linear change
of coordinates  $T$ of $ \proj{3}$,  the resulting parametrization
$T\circ \phi$ will satisfy B1 -- B5.  That is, if we allow generic
linear changes of coordinates of the image space, then B1 -- B4
still hold and B5 is a consequence of B1 -- B4.

Recall that $MP$ denotes both the map
\[\begin{CD}MP: R_{m-1,\,n-1}^4 @>{[a_0\; a_1\; a_2\; a_3]}>>
R_{2m-1,\,2n-1}\end{CD}\] given by
\[ (A_0,\, A_1,\, A_2,\, A_3) \mapsto \sum_{i=0}^3A_ia_i,\]
and the $4mn \times 4mn$ matrix which represents this map in the
standard monomial bases on $R_{m-1,\,n-1}^4 $ and
$R_{2m-1,\,2n-1}$. If we replace $\{a_i\}_{i=0}^3$ by
$\{a'_i\}_{i=0}^3$ where each $a'_i$ is a generic linear
combinations of $\{a_i\}_{i=0}^3$, then the rank of the
coefficient matrix $MP$ will not change.  Thus, the number of
linearly independent moving planes is also not affected by a
generic linear change of coordinates in the image space
$\proj{3}$.
\end{rem}
Let
\[ \begin{CD} MC: R_{m-1,\,n-1}^3 @>{[a_0\; a_1\; a_2]}>>
R_{2m-1,2n-1}
\end{CD}\] be the map given by
\[(A_0,\, A_1,\, A_2) \mapsto \sum_{i=0}^2A_ia_i.\] $MC$ is
represented by a matrix, also denoted $MC$ of size $4mn \times
3mn$, and $\Ker{MC}=\Syz(a_0,\, a_1,\, a_2)$. Thus
\[ \dim_{\C} \Syz(a_0,\, a_1,\, a_2)_{m-1,\,n-1}=
\dim_{\C} \Ker{MC}, \]
and the following fact is clear.
\begin{lem}
If $\phi: \proj{1} \times \proj{1} \to \proj{3}$  then $MC$ has
maximal rank ($=3mn$) if and only if $\phi$ satisfies condition
B6.
\end{lem}

We start our analysis with the following lemma, which indicates
that  base points of total multiplicity $k$ produce exactly $k$
linearly independent moving planes of bidegree $(m-1,\,n-1)$.
\begin{lem}\label{L:mp}
If  $\phi: \proj{1}\times \proj{1} \to \proj{3}$ satisfies the
base point conditions B1 -- B4, then
\[\dim_{\C} \Syz(I)_{m-1,\,n-1}=k.\]
\end{lem}

\begin{proof}
Consider the following exact sequence:
\[\begin{CD} 0 \rightarrow \Syz(I)_{m-1,\,n-1} \rightarrow
R^4_{m-1,\,n-1} @> {[a_0\;a_1\;a_2\;a_2\;a_3]}>> R_{2m-1,\,2n-1}
\rightarrow (R/I)_{2m-1,\,2n-1}\rightarrow 0. \end{CD}\] We have
\begin{eqnarray*}
\dim_{\C} \Syz(I)_{m-1,\,n-1} &=& \dim_{\C}(R/I)_{2m-1,\,2n-1}-
\dim_{\C}
R_{2m-1,\,2n-1} + 4 \dim_{\C} R_{m-1,\,n-1} \\
&=& \dim_{\C}(R/I)_{2m-1,\,2n-1}.
\end{eqnarray*}
Since at each base point, $\V(I)$ is a local complete
intersection, we have $\sum_{\mathbf{p}\in \V(I)} e(I,\mathbf{p})
= \deg(\V(I))=k$. Thus,
$\dim_{\C}(R/I)_{2m-1,\,2n-1}=\deg(\V(I))=k$, and hence
$$\dim_{\C} \Syz(I)_{m-1,\,n-1}=k.$$
\end{proof}

\begin{rem}
Under the hypotheses of Lemma \ref{L:mp}, the condition $\dim_{\C}
\Syz(I)_{m-1,\,n-1}=k$ means that there are exactly $k$ linearly
independent moving planes of bidegree $(m-1,\,n-1)$ which follow
the parametrization $\phi$.
\end{rem}

Our next goal is to prove that, under suitable conditions on the
base point scheme $\V(I)$,
\[\dim^{ }_{\C}  \Syz(I^2)_{m-1,\,n-1} =
mn+3k.\] We will start by proving the following two lemmas.
\begin{lem}\label{L:mqreg} If  $\phi:
\proj{1}\times \proj{1} \to \proj{3}$  satisfies the conditions B1
-- B4, and $I=\langle a_0,\, a_1,\, a_2,\, a_3 \rangle$ as usual,
then $\sat(I^2)$ is $(3m-1,\,3n-1)$-regular.

\end{lem}

\begin{proof} Consider the following exact sequence:
\[\begin{CD} 0 \rightarrow \Syz(I^2)_{m-1,\,n-1} \rightarrow
R^{10}_{m-1,\,n-1} @> {[a_0^2\; \cdots\; a_3^2]}>> R_{3m-1,\,3n-1}
\rightarrow (R/I^2)_{3m-1,\,3n-1}\rightarrow 0. \end{CD}\] This
implies that
\begin{eqnarray}\notag
\dim_{\C} \Syz(I^2)_{m-1,\,n-1} &=&
\dim_{\C}(R/I^2)_{3m-1,\,3n-1}- \dim_{\C} R_{3m-1,\,3n-1} + 10
\dim_{\C} R_{m-1,\,n-1} \\ \label{E:SyzII}&=&
\dim_{\C}(R/I^2)_{3m-1,\,3n-1}+mn.
\end{eqnarray}
Conditions B2, B3, and B4 show that
$\dim_{\C}(R/I)_{2m-1,\,2n-1}=\deg(\V(I))=k$, and this implies
that $I$ is $(2m-1,\,2n-1)$-regular by Theorem \ref{T:wkreg}.
Since $\V(I)$ is finite, Theorem \ref{T:satpower} shows that
$\sat(I^2)$ is
$((2-1)(2m-1)+m,\,(2-1)(2n-1)+n)=(3m-1,\,3n-1)$-regular, as
claimed.
\end{proof}

For the second lemma, we will need the following result of Herzog
\cite[Folgerung 2.2 and 2.4]{HZ}:

\begin{pro}\label{P:herzog} Let $\mc{O}_{\mbf{p}}$ be the
local ring of a point $\mbf{p}\in \proj{1}\times\proj{1}$,  and
let $\mc{I}_{\mbf{p}}\subseteq \mc{O}_{\mbf{p}}$ be a codimension
two ideal. Then
\begin{equation}\label{E:herzog}
\dim_{\mathbb{C}}\mc{I}_{\mbf{p}}/\mc{I}^2_{\mbf{p}} \ge
2\dim_{\mathbb{C}}\mc{O}_{\mbf{p}}/\mc{I}_{\mbf{p}},
\end{equation} and equality holds if and only if
$\mc{I}_{\mbf{p}}$ is a complete intersection ideal in
$\mc{O}_{\mbf{p}}$.
\end{pro}

\begin{lem}\label{L:mq}
If  $\phi: \proj{1}\times \proj{1} \to \proj{3}$  satisfies the
conditions B1 -- B4, then \[\dim_{\C}  \Syz(I^2)_{m-1,\,n-1} \geq
mn+3k.\]
\end{lem}

\begin{proof} The exact sequence
\[0 \rightarrow (I/I^2)_{r,\,r'}
\rightarrow (R/I^2)_{r,\,r'} \rightarrow (R/I)_{r,\,r'}
\rightarrow 0 \] shows that $\dim_{\C} (R/I^2)_{r,\,r'} =
\dim_{\C} (R/I)_{r,\,r'} + \dim_{\C} (I/I^2)_{r,\,r'},$  for all
$r$, $r'$. By condition B4 $\dim_{\C}
(R/I)_{2m-1,\,2n-1}=\deg(\V(I))=k,$ and since
$\dim_{\C}(R/I)_{k,\,k'}\le \dim_{\C}(R/I)_{l,\,l'}$ whenever
$(k,\,k')\le (l,\,l')$, it follows that $\dim_{\C}
(R/I)_{r,\,r'}=k$ for $r \geq 2m-1$, $r'\geq 2n-1$.  Hence,
$$\dim_{\C} (R/I^2)_{r,\,r'}=k+\dim_{\C}(I/I^2)_{r,\,r'}
\quad\text{for $r \geq 2m-1,\, r'\geq 2n-1$.}$$ For $r,r' \gg 0$,
$\dim_{\C} (I/I^2)_{r,\,r'} =P_{I/I^2}(r,\,r')$ where
$P_{I/I^2}(r,\,r')$ is the bigraded Hilbert polynomial of $I/I^2$.

If $\mc{I}$ is the sheaf of ideals associated to $I$, then
$\mc{I}/\mc{I}^2$ has  zero dimensional support since $\V(I)$ is
finite. Therefore, letting $X=\proj{1} \times \proj{1}$,
$$H^0(X,\,\mc{I}/\mc{I}^2)=
\bigoplus_{\mbf{p}\in
\V(I)}\mc{I}_{\mbf{p}}/\mc{I}_{\mbf{p}}^2\quad\text{and}\quad
H^0(X,\,\mc{I}/\mc{I}^2(r,\,r'))= \bigoplus_{\mbf{p}\in
\V(I)}\left(\mc{I}_{\mbf{p}}/\mc{I}_{\mbf{p}}^2\right)\otimes
\mc{O}_{\mbf{p}}(r,\,r'),$$ for all $r$, $r'$  and hence
\begin{equation*}
\dim_{\C} H^0(X,\,\mc{I}/\mc{I}^2)=\dim_{\C} H^0(X,\,
\mc{I}/\mc{I}^2(r,\,r'))\quad\text{for all $r$, $r'$,}
\end{equation*}
while $H^0(X,\, \mc{I}/\mc{I}^2(r,\,r'))=(I/I^2)_{r,\,r'}$ for all
$r,r' \gg 0$ by \cite[Theorem 1.6]{HY}. Therefore, for all $r,r'
\gg 0$ we have
\begin{eqnarray*}
P_{I/I^2}(r,\,r')&=&\dim_{\C}(I/I^2)_{r,\,r'}\\ &=&\dim_{\C}
H^0(X,\, \mc{I}/\mc{I}^2(r,\,r'))\\ &=&\dim_{\C}H^0(X,\,
\mc{I}/\mc{I}^2)\\ &=& \sum_{\mbf{p} \in \V(I)} \dim_{\C}
\mc{I}_{\mbf{p}}/ \mc{I}^2_{\mbf{p}}.
\end{eqnarray*}
Since each base point $\mbf{p}\in \V(I)$ is a local complete
intersection by condition B3, Proposition \ref{P:herzog} shows
that
\begin{equation}
 \sum_{\mbf{p}
\in \V(I)} \dim _{\C} \mc{I}_{\mbf{p}}/ \mc{I}^2_{\mbf{p}} = 2
\sum_{\mbf{p} \in \V(I)} \dim _{\C}
\mc{O}_{\mbf{p}}/\mc{I}_{\mbf{p}}=2\deg(\V(I))=2k,
\end{equation}
and hence,  for $r,\,r'\gg 0$,
\begin{equation}\label{E:3k}
\dim_{\C} (R/I^2)_{r,\,r'} = \dim_{\C} (R/I)_{r,\,r'} + \dim_{\C}
(I/I^2)_{r,\,r'} =k+ 2 \sum_{p \in \V(I)} \dim _{\C}
\mc{O}_p/\mc{I}_p = 3k.
\end{equation}

Since $P_{M^{\sat}}(r,\,r')=P_{M}(r,\,r')$ for any finitely
generated bihomogeneous $R$-module $M$, it follows  that
$$\dim_{\C} (R/I^2)_{r,\,r'} = \dim_{\C} (R/\sat(I^2))_{r,\,r'},$$
for $r,\,r' \gg 0$. This fact, combined with Equation
\eqref{E:3k}, the fact that $\sat(I^2)$ is $(3m-1,\,3n-1)$-regular
(Lemma \ref{L:mqreg}), and Lemma \ref{L:wkregdim0} shows that
$$\dim_{\C} (R/\sat(I^2))_{3m-1,\,3n-1}=3k.$$

Since $I^2 \subset \sat(I^2)$, we have $\dim_{\C}
(R/I^2)_{3m-1,\,3n-1} \geq \dim_{\C}
(R/\sat(I^2))_{3m-1,\,3n-1}=3k$. Therefore, Equation
\eqref{E:SyzII} becomes
\[\dim_{\C}  \Syz(I^2)_{m-1,\,n-1} = mn+
\dim_{\C}(R/I^2)_{3m-1,\,3n-1} \geq mn+3k. \]
\end{proof}

\begin{rem}
Under the hypothesis of Lemma \ref{L:mq}, the condition
\[\dim_{\C}\Syz(I^2)_{m-1,\,n-1}\geq mn+3k\] means that
there are at least $mn+3k$ linearly independent moving quadrics of
bidegree $(m-1,\,n-1)$ which follow the parametrization $\phi$.
\end{rem}

The construction of the matrix $M$ whose determinant is the
implicit equation of $S=\im{\phi}$ requires a careful choice of
basis of the vector space of moving quadrics, which is facilitated
by the following elementary linear algebra lemma.  We will first
establish the notation.

Let the vector space $V=V_1\oplus V_2$ be the direct sum of two
subspaces $V_1$ and $V_2$, and let $W \subset V$ be a subspace
such that $V_1 \cap W=\{0\}$.  Then the  projection $\pi: V \to
V_2$ along $V_1$ satisfies $\Ker \pi =V_1$, and $\Ker \pi |_{W} =
W \cap V_1=\{0\}$.  In particular, $\pi |_{W}$ is injective, so
that $\dim_{\C} W=\dim_{\C} \pi(W):=k$.  Let $\mc{B} =\{ v_1,
\ldots,\, v_l\}$ be a given basis of $V_2$.

\begin{lem}\label{L:linalg}
There is a subset $\mc{B}_1=\{v_{h_1},\, \ldots,\, v_{h_k}\}
\subset \mc{B}$ and a basis $\mc{C}=\{w_1,\, \ldots,\, w_k \}$ of
$W$ such that
\[\pi(w_e) =v_{h_e}+\overline{w}_e, \  \ \text{where }
\overline{w}_e \in
\Span (\mc{B} \setminus \mc{B}_1).\]
\end{lem}

\begin{proof}
Let $\{\widetilde{w}_1,\, \ldots,\, \widetilde{w}_k \}$ be an
arbitrary basis of $W$.   Then
\[ \pi( \widetilde{w}_i)=\sum^l_{j=1} a_{ij}v_j.\]
Let $A=[a_{ij}]$.  Then multiply $A$ on the left by an invertible
matrix  $P$ so that $PA=Q$, where $Q$ is in reduced row echelon
form.  Since $A$ is a $k \times l$ matrix which has $\rank A=k$
(because $\dim_{\C} \pi(W)=k$), there are $k$ columns $h_1 < h_2 <
\cdots <h_k$ which contain a leading $1$ in rows 1 to $k$,
respectively. Let $\mc{B}_1 = \{ v_{h_1},\,  \ldots,\, v_{h_k}\}$.
Let the basis $\mc{C}=\{w_1,\, \ldots,\, w_k\}$ be defined by
\[ \left[
\begin{matrix} w_1 \\ \vdots \\ w_k \end{matrix} \right] = P
\left[
\begin{matrix} \widetilde{w}_1 \\ \vdots \\ \widetilde{w}_k \end{matrix}
\right] ,                                    \] i.e., $w_e =\sum
^k_{j=1} p_{ej}\widetilde{w}_j$.  Then
\begin{eqnarray*}
\pi(w_e) & =& \sum^k_{j=1}p_{ej}\pi(\widetilde{w}_j) \\
&=& \sum^k_{j=1}p_{ej}  \sum^l_{r=1}a_{jr}v_r\\
& =& \text{Row}_e Q \left[
\begin{matrix} v_1 \\ \vdots \\ v_l \end{matrix} \right] \\
&=& v_{h_e} + \overline{w}_e
\end{eqnarray*}
where $\overline{w}_e \in \Span(\mc{B} \setminus  \mc{B}_1)$.
\end{proof}

If $P=\sum_{i=0}^3 A_i(s,\,u,\,t,\,v)x_i \in
R[x_1,\,x_2,\,x_3,\,x_4]$  is any moving plane, and
$L(x_0,\,x_1,\,x_2,\,x_3)$ is any homogeneous linear polynomial.
Then $P \cdot L$ is a moving quadric. Moreover, if $P$ follows
$\phi$, then $P\cdot L$ also follows $\phi$.  If $\mc{P}$ is a set
of moving planes, then $\mc{P} \cdot L :=\{P\cdot L : P \in \mc{P}
\}$. Let $\mc{P}_{\phi,\, m-1,\,n-1}$ be the set of moving planes
of bidegree $m-1,\,n-1$ which follow $\phi$, i.e., $(A_0,\, A_1,\,
A_2,\, A_3)_{m-1,\,n-1} \in \Syz(a_0,\, a_1,\, a_2,\,
a_3)_{m-1,\,n-1}$.

\begin{lem}\label{L:basismp}
Let $\phi:\proj{1}\times\proj{1}\to \proj{3}$, and assume that
$\phi$ satisfies condition B6, so that $\Syz(a_0,\, a_1,\,
a_2)_{m-1,\,n-1}=\{0\}$. Let $\mc{S}=\mc{P}_{\phi,\,m-1,\,n-1}$,
and let $\dim_{\C} \mc{S} =k$. Then $\mc{Q}=\sum_{i=0}^3
\mc{S}x_i$ is a vector space of moving quadrics which follow
$\phi$, with $\dim_{\C} \mc{Q}=4k$.
\end{lem}

\begin{proof}
We will apply Lemma \ref{L:linalg} with  the following
identifications:
\begin{itemize}
\item $V=\sum_{i=0}^3(R_{m-1,\,n-1})x_i\cong R^4_{m-1,\,n-1}$,
\item $V_1 =\sum_{i=0}^2(R_{m-1,\,n-1})x_i\cong R^3_{m-1,\,n-1}$,
\item
$V_2=(R_{m-1,\,n-1})x_3\cong R_{m-1,\,n-1}$, \item $W=\mc{S}$, and
\item $\mc{S} \cap V_1=\Syz(a_0,\, a_1,\, a_2)_{m-1,\,n-1} = \{0\}$.
\end{itemize}
Let $\mc{B}=\{ s^{\alpha} t^{\beta}x_3: 0 \leq \alpha \leq m-1,\,
0 \leq \beta \leq n-1\}$. According to Lemma \ref{L:linalg}, there
is a set $B=\{(\alpha_1,\,\beta_1),\, \ldots,\,
(\alpha_k,\,\beta_k)\}$ and a basis $\mc{C}=\{P_1,\, \ldots,\,
P_k\}$ of $\mc{S}$ such that $\pi(P_i)$, which is the part of
$P_i$ in $(R_{m-1,\,n-1})x_3$, has the form
\[\pi(P_i)= s^{\alpha_i} t^{\beta_i} x_3 +
\sum_{(\alpha,\,\beta)\notin B}b_{i,\,\alpha\,\beta} s^{\alpha}
t^{\beta}x_3\] where $i=1,\, 2,\, \ldots,\, k$. We claim that
$\{P_i x_j\}^{i=k,\,j=3}_{i=1,\,j=0}$ is a linearly independent
set. We  need to show that if
\begin{equation}\label{E:linquad}
\sum_{i=1}^k
\sum_{j=0}^3c_{ij} P_ix_j=0
\end{equation}
where $c_{ij} \in \C$, then we must have $c_{ij}=0$ for all $i$,
$j$. Since $$\mc{B}''=\{ s^{\alpha} t^{\beta}x_ix_j: 0 \leq \alpha
\leq m-1,\, 0 \leq \beta \leq n-1,\, 1\le i\le j\le 3\}$$ is a
basis of $ \bigoplus_{0\le i\le j\le 3}(R_{m-1,\,n-1})x_ix_j$, and
since $P_i$ is the only element of $\mc{C}$ that contains the term
$s^{\alpha_i}t^{\beta_i}x_3$, it follows that $P_ix_j$ is the only
term in \eqref{E:linquad} that contains the basis element
$s^{\alpha_i}t^{\beta_i}x_jx_3$ and hence the coefficient of this
term, namely $c_{ij}$, must be 0. Thus $c_{ij}=0$ for $i=1$, $2$,
$\ldots$, $ k$ and $j=0$, $1$, $2$, and $3$. Therefore, the moving
quadrics coming from the moving planes that follow $\phi$ are
linearly independent, and hence $\dim_{\C} \mc{Q}=4k$.
\end{proof}

\begin{thm}\label{T:rankMQ}
Let $ \phi: \proj{1} \times \proj{1} \to \proj{3}$ be given by
$\phi(s,u; t,v) =[a_0,\, a_1,\, a_2,\,a_3]$ where each $a_i\in R$
is a bihomogeneous polynomial of bidegree $(m,n)$,  and assume
that the base point scheme of $\phi$ satisfies conditions B1 - B6.
Then
$$\dim_{\C} \Syz(I^2)_{m-1,\,n-1}=mn+3k.$$
\end{thm}

\begin{proof}
If $MQ: R^{10}_{m-1,\,n-1} \to R_{3m-1,\,3n-1}$ is the map such
that $MQ(A_{00},\,A_{01},\,\ldots,\,A_{33}) = \sum_{0\leq i \leq j
\leq 3}A_{ij}a_ia_j$, we have that
\[10mn-\rank (MQ) =\dim_{\C} \Syz(I^2)_{m-1,\,n-1}\] is
the number of linearly independent moving quadrics. If $\rank (MQ)
\geq 9mn-3k$, then $$\dim_{\C} \Syz(I^2)_{m-1,\,n-1}=10mn-\rank
(MQ)\le mn+3k.$$ But  Lemma \ref{L:mq} shows that $\dim_{\C}
\Syz(I^2)_{m-1,\,n-1}\geq mn+3k$, and hence $\dim_{\C}
\Syz(I^2)_{m-1,\,n-1}= mn+3k$ will follow, once we have shown that
$\rank (MQ) \geq 9mn-3k$.

We now verify that this inequality is valid.   Since $\phi$
satisfies condition B6, the proof of Lemma \ref{L:basismp}, shows
that there is an indexed set $B=\{(\alpha_1,\, \beta_1),\,
\ldots,\, (\alpha_k,\, \beta_k)\}$ and a basis of moving planes
$\{P_1,\, \ldots,\, P_k\}$ such that
\begin{equation}\label{E:kmp}
P_i= s^{\alpha_i} t^{\beta_i} x_3 + \sum_{(\alpha,\,\beta)\notin
B}b_{i,\,\alpha\,\beta} s^{\alpha} t^{\beta}x_3
+\sum_{j=0}^2\sum_{(\alpha,\,\beta)}c_{i,\,\alpha\,
\beta}s^{\alpha}t^{\beta}x_j
\end{equation}
where $i=1$, $ 2$, $\ldots$, $ k$.  As with $MP$, the matrix
representing $MQ$ with respect to the standard bases is also
denoted $MQ$. Thus  the columns of $MQ$ are indexed by
\[\Lambda=\{ s^{\alpha}t^{\beta}x_ix_j: 0 \leq \alpha \leq m-1,\,
0 \leq \beta \leq n-1 ,\, 0\leq i \leq j \leq 3 \}.\] If
\[
  \Lambda_P=\{
  s^{\alpha_i}t^{\beta_i}x_jx_3,\, s^{\alpha}t^{\beta}x_3^2
:  1 \leq i \leq k,\, 0 \leq \alpha \leq m-1,\, 0 \leq \beta \leq
n-1,\, 0 \leq j \leq 2 \},
\]
and $\Lambda'=\Lambda \setminus \Lambda_P$, then
$|\Lambda'|=10mn-(mn+3k)=9mn-3k$.   Let $MQ'$ be the matrix
obtained from $MQ$ by deleting the columns indexed by $\Lambda_p$.
Thus the nonzero  elements of $\Ker{MQ'}$ correspond to nontrivial
syzygies:
\begin{equation}\label{E:MQ'}
A_{00}a_0^2+A_{01}a_0a_1+A_{02}a_0a_2+A_{03}a_0a_3+A_{11}a_1^2+
A_{12}a_1a_2+A_{13}a_1a_3 +A_{22}a_2^2+A_{23}a_2a_3=0
\end{equation} where $A_{ij}$ is bihomogeneous of bidegree
$(m-1,\,n-1)$ and there are no terms $ s^{\alpha_i}t^{\beta_i}$ in
$\{A_{i3}\}_{i=0}^{2}$. Since every term contains $a_0$,  $a_1$,
or $a_2$, we obtain:
\[(A_{00}a_0+A_{01}a_1+A_{02}a_2+A_{03}a_3)a_0+(A_{11}a_1+A_{12}
a_2+A_{13}a_3)a_1+ (A_{22}a_2+A_{23}a_3)a_2=0,\] which means that
$$(B_1,\, B_2,\,B_3)=(A_{00}a_0+A_{01}a_1+A_{02}a_2+A_{03}a_3,\,
A_{11}a_1+A_{12}a_2+A_{13}a_3,\, A_{22}a_2+A_{23}a_3)$$ is a
syzygy of $\langle a_0,\, a_1,\, a_2\rangle$.  Each $B_i$ has
bidegree $(2m-1,\,2n-1)$ and
$$B_i\in \langle a_0,\,a_1,\,a_2,\,a_3\rangle \subseteq
\sat\langle a_1, a_1,\,a_2\rangle$$ by condition B5.  Therefore,
each $B_i$ vanishes on the base point scheme of $\langle a_0,\,
a_1,\, a_2\rangle$. By condition B5, $\V(a_0,\,a_1,\,a_2)=\V(I)$,
and thus, by B2,  $\V( a_0,\,a_1,\,a_2)\subset \proj{1}\times
\proj{1}$ is finite and each base point is a local complete
intersections (condition B3). Since each $a_i$ has bidegree
$(m,\,n)$, while each $B_i$ has bidegree $(2m-1,\,2n-1)$, we have
$(2m-1 -2m+1)(2n-1-2n+1)= 0.$ Therefore, all the hypotheses of
Theorem \ref{T:lcikoz} are satisfied, and we conclude that all
syzygies of bidegree $(2m-1,\,2n-1)$ of $\langle
a_0,\,a_1,\,a_2\rangle$ that vanish on the base point scheme are
in fact Koszul syzygies. Since $(B_0,\,B_1,\,B_2)$ is a syzygy
that vanishes on $\V(a_0,\,a_1,\,a_2)$, it follows that there are
bihomogeneous polynomials $h_1$, $h_2$, and $h_3$ in $R$ of
bidegree $(m-1,\,n-1)$ such that:
\begin{eqnarray*}
A_{00}a_0+A_{01}a_1+A_{02}a_2+A_{03}a_3&=&h_1a_2+h_2a_1 \\
A_{11}a_1+A_{12}a_2+A_{13}a_3&=&-h_2a_0+h_3a_2 \\
A_{22}a_2+A_{23}a_3&=&-h_1a_0-h_3a_1.
\end{eqnarray*}
We can rewrite the above equations to get:
\begin{eqnarray}\label{E:3syz}
 A_{00}a_0+(A_{01}-h_2)a_1+(A_{02}-h_1)a_2+A_{03}a_3&=&0, \\
\label{E:3syz1}
 h_2a_0+A_{11}a_1+(A_{12}-h_3)a_2+A_{13}a_3&=&0,\\ \label{E:3syz2}
 h_1a_0+h_3a_1+A_{22}a_2+A_{23}a_3&=&0.
\end{eqnarray}
We know that $A_{ij}$ is bihomogeneous of bidegree $(m-1,\,n-1)$
and there are no $s^{\alpha_i}t^{\beta_i}$ terms in
$\{A_{i3}\}_{i=0}^{2}$. Thus Equations \eqref{E:3syz},
\eqref{E:3syz1}, \eqref{E:3syz2} are nontrivial  syzygies of
$\langle a_0,\, a_1,\, a_2,\, a_3\rangle$ which correspond to
moving planes $P$ with no $s^{\alpha_i}t^{\beta_i}x_3$ term for $1
\leq i \leq k$. But $\{P_1,\, \ldots,\, P_k\}$ is a basis of
moving planes. Any nonzero moving plane $P=c_1P_1 +\cdots +
c_kP_k$ must have some nonzero term $s^{\alpha_i}t^{\beta_i}x_3$,
since if $c_i \neq 0$, then $s^{\alpha_i}t^{\beta_i}x_3$ appears.

Hence the nontrivial syzygies from Equations \eqref{E:3syz},
\eqref{E:3syz1}, \eqref{E:3syz2} cannot exist. Thus
$\Ker{MQ'}=\{0\}$, so $$\rank( MQ) \geq \rank( MQ') =9mn-3k,$$ as
required, and hence we conclude that $\dim_{\C}
\Syz(I^2)_{m-1,\,n-1}=mn+3k$.
\end{proof}

\begin{rem}
Under the hypothesis of Theorem \ref{T:rankMQ}, the condition $
\Syz(I^2)_{m-1,n-1} = mn+3k$ means that there are exactly $mn+3k$
linearly independent moving quadrics of bidegree $(m-1,\,n-1)$
that follow the parametrization $\phi$. Moreover, the proof shows
that there are no nontrivial moving quadrics with nonzero
coordinates coming only from the basis elements
$\Lambda'=\Lambda\setminus\Lambda_P$ (because $\Ker{MQ'}=\{0\}$).
Hence any nontrivial moving quadric $Q$ must have at least one
nonzero coordinate from a term in the set
\[
  \Lambda_P=\{
  s^{\alpha_i}t^{\beta_i}x_jx_3,\, s^{\alpha}t^{\beta}x_3^2
:  1 \leq i \leq k,\, 0 \leq \alpha \leq m-1,\, 0 \leq \beta \leq
n-1,\, 0 \leq j \leq 2 \}.
\] This observation will be key to the proof of Theorem
\ref{T:surimp}.
\end{rem}

If $\phi:\proj{1}\times\proj{1}\to \proj{3}$ satisfies the base
point conditions B1 -- B6, then Lemma \ref{L:mp} shows that there
are exactly $k$ linearly independent moving planes
$\mc{MP}=\{P_\gamma\}_{\gamma=1}^k$ which follow $\phi$, where
$k\le mn$ is the total multiplicity of all base points of $\phi$,
and Theorem \ref{T:rankMQ} shows that there are exactly $mn+3k$
linearly independent moving quadrics $\mc{MQ}=\{Q_{\tau}:1\le
\tau\le mn+3k\}$ that follow $\phi$.  Each moving plane can be
written as
$$P_{\gamma}=\sum_{i=0}^3A_ix_i=\sum_{\alpha=0}^{m-1}\sum_{\beta=0}^
{n-1}P_{\gamma,\alpha\,\beta}(x_0,\,x_1,\,x_2,\,x_3)s^{\alpha}t^{\beta},
$$
and each  moving quadric $Q_\tau$  can be written as (see Equation
\eqref{E:quad})
\begin{equation*}
Q_\tau =  \sum_{0\leq i \leq j \leq 3} A_{ij} x_ix_j=
\sum_{\alpha=0}^{m-1}\sum_{\beta=0}^{n-1}Q_{\tau,\,\alpha\, \beta}
(x_0,x_1,x_2,x_3)s^{\alpha}t^{\beta},
\end{equation*} where $P_{\gamma,\alpha\,\beta}(x_0,\,x_1,\,x_2,\,x_3)$
is a homogeneous linear form and $Q_{\tau,\,\alpha\, \beta}$ is a
homogeneous quadratic form in $x_i$ with coefficients in $\C$. Our
goal is to choose the sets of moving planes $\mc{MP}$ and moving
quadrics $\mc{MQ}$ in such a way that all $k$ of the moving planes
and $mn-k$ of the moving quadrics can be combined into a single
$mn\times mn$ matrix (which will depend on the choice of $\mc{MP}$
and $\mc{MQ}$)
\begin{equation}\label{E:MPQ}
M= \begin{bmatrix}
  P_{\gamma,\,\alpha\,\beta}(x_0,\,x_1,\,x_2,\,x_3) \\
  Q_{\tau,\,\alpha\,\beta}(x_0,\,x_1,\,x_2,\,x_3)
\end{bmatrix}
\end{equation}
such that the equation of the image surface $S=\im{\phi}$ is given
by the determinantal equation $|M|=0$, as long as $\phi$ is
generically one-to-one.  The strategy for constructing $M$ is to
start with an arbitrary basis $\mc{MP}$ of moving planes
(consisting of $k$ moving planes), and then choose a basis of
moving quadrics $\mc{MQ}$ (consisting of $mn+3k$ moving quadrics)
in such a manner that $4k$ of the moving quadrics are obtained by
multiplying the moving planes of $\mc{MP}$ by each of the
coordinate functions $x_i$ ($0\le i\le 3$).  If these $4k$ moving
quadrics are deleted from the set $\mc{MQ}$, then the remaining
$mn-k$ are used for the matrix $M$ of Equation \eqref{E:MPQ}. The
justification for this procedure constitutes the proof of our main
result.

\begin{thm}\label{T:surimp}
Let $\phi: \proj{1}\times \proj{1} \to \proj{3}$ be a
parametrization of a surface $S=\im{\phi}\subset \proj{3}$. If
$\phi$ is generically one-to-one and satisfies base point
conditions B1 -- B6, then a basis of moving planes $\mc{MP}$ and a
basis of moving quadrics $\mc{MQ}$ can be chosen so that $S$ is
defined by the determinantal equation $|M|=0$, where $M$ is the
$mn\times mn$ matrix of Equation \eqref{E:MPQ}.
\end{thm}

\begin{proof} By Lemma \ref{L:mp}, $\dim_{\C} \Syz(a_0,\, a_1,\, a_2,
\,a_3)_{m-1,n-1}=k$, and the proof of Lemma \ref{L:basismp} shows
that there is an indexed set $B=\{(\alpha_1,\, \beta_1),\,
\ldots,\, (\alpha_k,\, \beta_k)\}$ and a basis of moving planes
$\mc{MP}=\{P_1,\, \ldots,\, P_k\}$ such that
\begin{equation}\label{E:kmp2}
P_i= s^{\alpha_i} t^{\beta_i} x_3 + \sum_{(\alpha,\,\beta)\notin
B}b_{i,\,\alpha\,\beta} s^{\alpha} t^{\beta}x_3
+\sum_{j=0}^2\sum_{(\alpha,\,\beta)}c_{i,\,\alpha\,
\beta}s^{\alpha}t^{\beta}x_j
\end{equation}
where $i=1$, $ 2$, $\ldots$, $ k$.   By Theorem \ref{T:rankMQ},
$$\dim_{\C} \Syz(I^2)_{m-1,\,n-1}=mn+3k.$$
We now describe how  to produce a convenient basis of
$\Syz(I^2)_{m-1,\,n-1}$.

Let $V_{\Lambda'}$,$V_{\Lambda_p}$ be subspaces of
$V=\bigoplus_{0\le i \le j\le
3}\left(R_{m-1,\,n-1}\right)x_ix_j\cong R^{10}_{m-1,\,n-1}$ with
bases $\Lambda'$, $\Lambda_P$ respectively, where (as in the proof
of Theorem \ref{T:rankMQ})
\[\Lambda=\{ s^{\alpha}t^{\beta}x_ix_j: 0 \leq \alpha \leq m-1,\;
0 \leq \beta \leq n-1,\; 0 \leq i \leq j \leq 3 \},\]
\[  \Lambda_P=\{ s^{\alpha_i}t^{\beta_i}x_jx_3,\;
s^{\alpha}t^{\beta}x_3^2:  0 \leq i \leq k,\; 0 \leq \alpha \leq
m-1,\; 0 \leq \beta \leq n-1,\; 0\leq j \leq 2 \}.
\]
and $\Lambda'=\Lambda \setminus \Lambda_P$. Then
$V=V_{\Lambda'}\oplus V_{\Lambda_p}$ and the proof of Theorem
\ref{T:rankMQ} (namely, the proof that $\Ker{MQ'}=\{0\}$) shows
that $\Syz(I^2)_{m-1,\,n-1} \subset V$ satisfies
$$\Syz(I^2)_{m-1,\,n-1} \cap V_{\Lambda'} =\{0\}.$$ We conclude
that if $\pi: V \to V_{\Lambda_P}$ given by $\pi(v_1+v_2)=v_2$ is
the projection onto $V_{\Lambda_P}$ along $V_{\Lambda'}$, then
$\pi|_{\Syz(I^2)_{m-1,n-1}}$ is an isomorphism, since $$\dim_{\C}
\Syz(I^2)_{m-1,n-1}=\dim_{\C} V_{\Lambda_P}=mn+3k.$$ Thus
$\mc{MQ}=\pi^{-1}(\Lambda_P)$ is a basis of moving quadrics that
follow $\phi$.

Let $Q_{x_jx_3,i}=\pi^{-1}(s^{\alpha_i}t^{\beta_i}x_jx_3)$, for $1
\leq i \leq k$, and $0 \leq j \leq 3$. Since $x_jP_i \in
 \Syz(I^2)_{m-1,n-1}$, and
$\pi(x_jP_i)=s^{\alpha_i}t^{\beta_i}x_jx_3$ (see Equation
\ref{E:kmp2}), the fact that $\pi|_{\Syz(I^2)_{m-1,n-1}}$ is an
isomorphism shows that  $x_jP_i=Q_{x_jx_3,i}$.  Thus, we have
identified the set of moving quadrics in $\mc{MQ}$ which arise
from multiplication of the moving planes in $\mc{MP}$ by the
homogeneous coordinate functions $x_j$ ($0\le j\le 3$).  These are
excluded when forming the matrix $M$.

Let $Q_{\gamma\,\delta}=\pi^{-1}(s^{\gamma}t^{\delta}x_3^2)$,
where $$( \gamma , \delta) \in \{(\alpha, \beta):0 \leq \alpha
\leq m-1,\; 0 \leq \beta \leq n-1\} \setminus \{(\alpha_i,
\beta_i): 1 \leq i \leq k \}:=C_P.$$  These $mn-k$ moving quadrics
in the basis $\mc{MQ}$ do not come from the moving planes of
$\mc{MP}$ by multiplication by $\{x_i\}_{i=0}^3$. Thus, they can
be combined with the $k$ moving planes $\mc{MP}$ to produce the
matrix $M$. Hence
\begin{equation}\label{E:MPQ2}
M=\begin{bmatrix}
  P_{i,\,\alpha\,\beta}(x_0,\,x_1,\,x_2,\,x_3) \\
  Q_{\gamma\,\delta,\,\alpha\,\beta}(x_0,\,x_1,\,x_2,\,x_3)
\end{bmatrix}
\end{equation}
where $1\le i\le k$ and $(\gamma,\,\delta)\in C_P$, and the
columns are indexed by the monomial basis $s^\alpha t^\beta$ of
$R_{m-1,\, n-1}$ with $0\le \alpha\le m-1$, $0\le \beta\le n-1$.

The $k$ moving planes $P_i$ have the form
\[P_i= s^{\alpha_i} t^{\beta_i} x_3 + \sum_{(\alpha,\,\beta)\notin
B}b_{i,\,\alpha\,\beta} s^{\alpha} t^{\beta}x_3
+\sum_{j=0}^2\sum_{(\alpha,\,\beta)}c_{i,\,\alpha\,
\beta}s^{\alpha}t^{\beta}x_j,
\] while the $mn-k$ moving quadrics $Q_{\gamma\,\delta}$ for
$(\gamma,\,\delta)\in C_P$ have the form
\[ Q_{\gamma \delta}=x_3^2s^{\gamma}t^{\delta}+
\text{ terms not involving }  x_3^2.\] That is, the term
$s^{\alpha_i} t^{\beta_i} x_3$ occurs in $P_i$, but in no other
$P_j$ for $j\ne i$, while the term $x_3^2s^{\gamma}t^{\delta}$
occurs in $ Q_{\gamma\, \delta}$, but  no other term of the form
$x_3^2s^{\gamma'}t^{\delta'}$ occurs in $Q_{\gamma\,\delta}$. Thus
the matrix $M$ of Equation \eqref{E:MPQ2} will have $k$ linear
rows and $mn-k$ quadratic rows in the variables $x_0$, $x_1$,
$x_2$ and $x_3$. Moreover, we can order the rows and columns in
such a way that all of the $x_3^2$ terms (one for each quadratic
row) occur on the last $mn-k$ diagonals, while the the first $k$
diagonals have the term $x_3$ coming from the terms $s^{\alpha_i}
t^{\beta_i} x_3 $ in $P_i$ ($1\le i\le k$).  Thus, after
appropriate ordering of the rows and columns, $M$ will have the
form

\[ {M}=\left[
\begin{matrix}
 x_3+\cdots & \\
    & \ddots \\
    &  & x_3+\cdots & \\
                 &             &          &
x_3^2+\cdots\\
               &             &          &
        &  \ddots\\
                &             & &  &&
x_3^2+\cdots
 \end{matrix}\right].
\]
There are $k$ linear rows and $mn-k$ quadratic rows, so the
determinant of $M$  contains the term $x_3^{2mn-k}$, which occurs
in the multiplication of the diagonal entries.  Since the $x_3^2$
term appears only in the last $mn-k$ diagonal entries, and in the
upper left $k\times k$ block, the term $x_3$ appears only on the
diagonal, it follows that $2mn-k$ is the highest power of $x_3$
that can appear in $|M|$, and this power appears with nonzero
coefficient. Thus $|{M}|$ is not identically zero. Since $M$
contains $mn-k$ rows of quadratic terms in $x_i$ and $k$ rows of
linear terms in $x_i$, the total degree of $|{M}|$ is $2mn-k$. By
construction, the rows of $M$ represent moving quadrics  and
moving planes that follow the surface, and hence, when $x_i$ is
replaced by $a_i$ it follows that the columns of $M$ are linearly
dependent. Therefore, $|M|$ vanishes for points on the surface.
From the degree formula
\[ \text{deg}(\phi) \text{
deg}(S)=2mn-\sum_{\text{base points}}\text{multiplicity of the
base point}
\]
and the fact that $\phi$ is generically one to one, and the total
multiplicity  of all base points is $k$, we conclude that
$\text{deg}S=2mn-k$. But this is the same degree as $|M|$ so
$|M|=0$ must be the implicit equation of the image of $\phi$.
\end{proof}

\begin{exa}
Consider $\phi:\proj{1}\times\proj{1}\to \proj{3}$  given by the
following parametrization:
\[
\begin{matrix}
a_0=u^2tv+s^2tv,& a_1=u^2t^2+suv^2,& a_2=s^2v^2+s^2t^2,&
a_3=s^2tv.
\end{matrix}
\] Here $m=n=2$,  $\V(a_0, a_1, a_2, a_3)=(0:1;0:1)$ and
the multiplicity of the single base point is one.  This base point
is a local complete intersection and $\deg \V(I)=1$. Using
Singular \cite{GPS01}, one can verify that the base point
conditions B1 -- B6 are satisfied, since $\dim_{\C}(R/I)_{3,3}=1$,
$a_3 \in \sat(a_0, a_1,a_2)$ and $\dim_{\C}\Syz(a_0,
a_1,a_2)_{1,1}=0$. Also, by Singular, we find one moving plane of
bidegree $(1,1)$ which is
\[ -x_2+tx_3+sx_1+st(x_3-x_0)\] and three linearly independent
moving quadrics of bidegree $(1,1)$ which are complementary to
$x_i(-x_2+tx_3+sx_1+st(x_3-x_0))$  for $i=0,1,2,3$:
\[ x_0x_3+t(x_0x_2+x_1x_3+x_2x_3)+s(-x_0x_3+x_3^2) \]
\[ (x_1x_3-x_2x_3)+t(x_0x_3+2x_3^2)+s(-x_0x_2+x_1x_3+x_2x_3) \]
\[ (x_2^2-x_3^2)+t(-x^2x_3)+s(x_0x_3-x_1x_2-x_3^2)+stx_1x_3\]
Thus the matrix $M$ is
\[ M= \left[ \begin{matrix} -x_2&x_3&x_1&x_3-x_0\\
x_0x_3&x_0x_2+x_1x_3+x_2x_3&-x_0x_3+x_3^2&0\\
x_1x_3-x_2x_3&x_0x_3+2x_3^2&-x_0x_2+x_1x_3+x_2x_3&0\\
x_2^2-x_3^2&-x^2x_3&x_0x_3-x_1x_2-x_3^2&x_1x_3 \end{matrix}
\right]
\] and
\begin{eqnarray*}
|M|&=&-x_0^3x_2^4+x_0^2x_1^2x_2^2x_3+x_0^2x_2^4x_4+
x_0x_1^2x_2^2x_3^2+\\
& &2x_0x_1x_2^3x_3^2+x_0x_2^4x_3^2-x_0^4x_3^3-2x_0^2x_1^2x_3^3-
x_1^4x_3^3+\\
& &3x_0^2x_1x_2x_3^3-2x_1^3x_2x_3^3+3x_0^2x_2^2x_3^3-2x_1^2x_2^2x_3^3-
2x_1x_2^3x_3^3\\
& &-x_2^4x_3^3+5x_0^3x_3^4+x_0x_1^2x_3^4-7x_0x_1x_2x_3^4-
6x_0x_2^2x_3^4-9x_0^2x_3^5\\
& &+x_1^2x_3^5+4x_1x_2x_3^5+3x_2^2x_3^5+7x_0x_3^6-2x_3^7.
\end{eqnarray*}
Thus the theorem gives $|M|=0$ as the implicit equation of
$S=\im{\phi}$. Note $|M|$ is a polynomial of degree $7$ which is
the same as the degree of the parametrized surface.
\end{exa}


\begin{thebibliography}{10}

\bibitem{BRS}
M.~P. Brodmann and R.~Y. Sharp.
\newblock {\em Local Cohomology}.
\newblock Cambridge studies in advanced mathematics. Cambridge University
  Press, 1998.

\bibitem{BCD}
L.~Bus\'{e}, D.~Cox, and C.~D'Andrea.
\newblock Implicitization of surfaces in $\proj{3}$ in the presence of base
  points.
\newblock {\em Journal of Algebra and its Applications}.
\newblock to appear.

\bibitem{KC}
K.~A. Chandler.
\newblock Regularity of the powers of an ideal.
\newblock {\em Comm. Algebra}, 25:3773--3776, 1997.

\bibitem{C1}
D.~A. Cox.
\newblock Curves, surfaces, and syzygies.
\newblock {\em Contemporary Mathematics}, 286:1--20, 2001.

\bibitem{C}
D.~A. Cox.
\newblock Equations of parametric curves and surfaces via syzygies.
\newblock {\em Contemporary Mathematics}, 286:1--20, 2001.

\bibitem{CGZ}
D.~A. Cox, R.~N. Goldman, and M.~Zhang.
\newblock On the validity of implicitization by moving quadrics for rational
  surfaces with no base points.
\newblock {\em J. Symb. Comput.}, 29:419--440, 2000.

\bibitem{GPS01}
G.-M. Greuel, G.~Pfister, and H.~Sch\"onemann.
\newblock {\sc Singular} 2.0.
\newblock {A Computer Algebra System for Polynomial Computations}, Centre for
  Computer Algebra, University of Kaiserslautern, 2001.
\newblock {\tt http://www.singular.uni-kl.de}.

\bibitem{HZ}
J.~Herzog.
\newblock {Ein Cohen-Macauly-Kriterium mit Anwedungen auf den Konormalenmodul
  und den Differentialmodule}.
\newblock {\em Math. Z.}, 163:149--162, 1978.

\bibitem{HW}
J.~W. Hoffman and H.~Wang.
\newblock Castelnuovo-{M}umford regularity in biprojective spaces.
\newblock {\tt http://arxiv.org/abs/math.AG/0212033}.

\bibitem{HW2}
J.~W. Hoffman and H.~Wang.
\newblock Curvilinear base points, local complete interesction and koszul
  syzygies in biprojective spaces.
\newblock {\tt http://arxiv.org/abs/math.AG/0304118}.

\bibitem{HY}
E.~Hyry.
\newblock {The diagonal subring and the Cohen-Macaulay property of a
  multigraded ring}.
\newblock {\em Trans. Amer. Math. Soc.}, 351(6):2213--2232, 1999.

\bibitem{DM}
D.~Mumford.
\newblock {\em Lectures on curves on an algebraic surface}.
\newblock Princeton University Press, Princeton, New Jersey, 1966.

\bibitem{O}
A.~Ooishi.
\newblock Castelnuovo's regularity of graded rings and modules.
\newblock {\em Hiroshima Math. J.}, 12:627--644, 1982.

\bibitem{SW}
J.~H. Sampson and G.~Washnitzer.
\newblock A {K\"u}nneth formula for coherent algebraic sheaves.
\newblock {\em Illinois J. Math.}, 3:389--402, 1959.

\bibitem{Wang}
H.~Wang.
\newblock {\em Equations of Parametric Surfaces with Base Points via Syzygies}.
\newblock PhD thesis, Louisiana State University, 2003.

\end{thebibliography}
\end{document}